\documentclass[numreferences]{kluwer}
\usepackage[T1]{fontenc}
\usepackage{geometry}
\usepackage[dvips]{color}
\usepackage{amssymb}

\makeatletter
 \newtheorem{thm}{Theorem}[section]
 \newtheorem{thm*}{Theorem} 
 \newdisplay{defn}[thm]{Definition}
 \newtheorem{prop}[thm]{Proposition} 
 \newdisplay{rem}[thm]{Remark}
 \newtheorem{lem}[thm]{Lemma} 
 \newdisplay{note}[thm]{Note} 
 \newdisplay{case}{Case} 
 \newtheorem{cor}[thm]{Corollary} 

\newcommand\rest{\hbox{\raise.17em\hbox{$ |\kern-.2em$}\lower.23em\hbox{$-$}}}
\newcommand\avint{\hbox{\hbox{$\displaystyle \int$}\hbox{\kern-.9em{$-$}}}}
\newcommand\smavint{\hbox{\hbox{$\int$}\hbox{\kern-.75em{$-$}}}}
\newcommand{\R}{\mathbb{R}}

\makeatother
\begin{document}
\begin{frontmatter}
\end{frontmatter}
\begin{opening}
\title{Partial regularity of mass-minimizing rectifiable sections }
\runningtitle{Partial regularity of mass-minimizing rectifiable sections }

\author{David L. \surname{Johnson} \email{david.johnson@lehigh.edu}}
\author{Penelope \surname{Smith} \email{ps02@lehigh.edu}}


\institute{Lehigh University, Bethlehem, PA}

\begin{ao}\\
Department of Mathematics\\
Lehigh University\\
14 E. Packer Avenue\\
Bethlehem, Pennsylvania 18015-3174
\end{ao}


\begin{abstract}
Let $B$ be a fiber bundle with compact fiber $F$ over a compact
Riemannian $n$-manifold $M^{n}$. There is a natural Riemannian metric
on the total space $B$ consistent with the metric on $M$. With respect
to that metric, the volume of a rectifiable section $\sigma:M\rightarrow B$
is the mass of the image $\sigma(M)$ as a rectifiable $n$-current
in $B$. 
\begin{thm*}
For any homology class of sections of $B$, there is a mass-minimizing
rectifiable current $T$ representing that homology class which is
the graph of a $C^{1}$ section on an open dense subset of $M$.
\end{thm*}
\end{abstract}


\classification{2000 Mathematics Subject Classification}{49F20, 49F22, 49F10, 58A25, 53C42, 53C65}

\end{opening}
\begin{article}

\section*{Introduction }

The notion of the volume of a section of a fiber bundle over a manifold
$M$ was introduced by H. Gluck and W. Ziller, in the special case
of the unit tangent bundle $\pi:T_{1}(M)\rightarrow M$, where sections
are unit vector fields, or flows on $M$. The volume of a section
$\sigma$ is defined as the mass (Hausdorff $n$-dimensional measure)
of the image $\sigma(M)$. They were able to establish, by constructing
a calibration, that the tangents to the fibers of the standard Hopf
fibration $S^{3}\rightarrow S^{2}$ minimized volume among all sections
of the unit tangent bundle of the round $S^{3}$.

However, in general calibrations are not available, even for the unit
tangent bundles of higher-dimensional spheres. For a general bundle
$\pi:B\rightarrow M$ over a Riemannian $n$-manifold $M$, with compact
fiber $F$, there is a special class of rectifiable currents, called
\emph{rectifiable sections}, which includes all smooth sections and
which has the proper compactness properties to guarantee the existence
of volume-minimizing rectifiable sections in any homology class. Partial
regularity of volume-minimizing rectifiable sections in general is
the subject of this paper. 

The basic partial-regularity result established here is that a volume-minimizing
rectifiable section exists in any homology class of sections which
is a $C^{1}$ section over an open, dense subset of $M$. This does
not state that a dense subset of the section itself consists of regular
points. In fact, there are simple counter-examples of that statement.
Denseness of the set of points in $M$ over which the section is regular
is straightforward, but openness in $M$ requires some work.

Our approach to this problem begins with a penalty functional, composed
of the $n$-dimensional area integrand plus a parameter ($1/\epsilon$)
multiplied by a term measuring the deviation from a graph of a current
in the total space. Each penalty functional will have energy-minimizing
currents which are rectifiable currents in the total space, but which
are not necessarily rectifiable sections. As the penalty parameter
$\epsilon$ approaches $0$, the {}``bad\char`\"{} set of points
in the base over which the current is not a section will have small
measure, and outside a slightly larger set the current will be a $C^{1}$
graph. These penalty minimizers will converge to a rectifiable section
which will be a minimizer of the volume problem. 

Once fundamental monotonicity properties are established for this
limiting minimizer, the program to establish partial regularity of
energy-minimizing currents due to Bombieri in \cite{Bombieri} can
be applied, with significant modifications for the current situation,
to show that the limiting minimizer is sufficiently smooth on an open
dense set. 

The main theorem of this paper is the following:

\begin{thm}
\label{thm:main} Let $B$ be a fiber bundle with compact fiber $F$
over a compact Riemannian manifold $M$, endowed with the Sasaki metric
from a connection on $B$. For any homology class of sections of $B,$
there is a mass-minimizing rectifiable section $T$ representing that
homology class which is the graph of a $C^{1}$ section on an open
dense subset of $M$. 
\end{thm}

\section{Definitions}

Let $B$ be a Riemannian fiber bundle with compact fiber $F$ over
a Riemannian $n$-manifold $M$, with projection $\pi:B\rightarrow M$
a Riemannian submersion. $F$ is a $j$-dimensional compact Riemannian
manifold. Following \cite{MS}, $B$ embeds isometrically in a vector
bundle $\pi:E\rightarrow M$ of some rank $k\geq j$, which has a
smooth inner product $<\:,\:>$ on the fibers, compatible with the
Riemannian metric on $F$. The inner product defines a collection
of connections, called \textit{metric connections}, which are compatible
with the metric. Let a metric connection $\nabla$ be chosen. The
connection $\nabla$ defines a Riemannian metric on the total space
$E$ so that the projection $\pi:E\rightarrow M$ is a Riemannian
submersion and so that the fibers are totally geodesic and isometric
with the inner product space $E_{x}\cong\R^{k}$ \cite{Sasaki}, \cite{j}. 

We will be using multiindices $\alpha=(\alpha_{1},\dots,\alpha_{n-l})$,
$\alpha_{i}\in\{1,\dots,n\}$ with $\alpha_{1}<\cdots<\alpha_{n-l}$,
over the local base variables, and $\beta=(\beta_{1},\dots,\beta_{l})$,
$\beta_{j}\in\left\{ 1,\dots,k\right\} $ with $\beta_{1}<\cdots<\beta_{l}$,
over the local fiber variables (we will at times need to consider
the vector bundle fiber, as well as the compact fiber $F$; which
is considered will be clear by context). The range of pairs $(\alpha,\beta)$
is over all pairs satisfying $|\beta|+|\alpha|=n$, where $\left|\left(\alpha_{1},\dots,\alpha_{m}\right)\right|:=m$.
As a notational convenience, denote by $n$ the $n$-tuple $n:=(1,\ldots,n)$,
and denote the null $0$-tuple by $0$. 

In addition, we indicate by $A\ll B$ that $A$ is bounded above by
a constant times $B$, where that constant is independent of the variables
included in $A$ and $B$.

\begin{defn}
An $n$-dimensional current $T$ on a Riemannian fiber bundle $B$
over a Riemannian $n$-manifold $M$ locally, over a coordinate neighborhood
$\Omega$ on $M$, decomposes into a collection, called \emph{components},
or \emph{component currents} \emph{of} $T$, with respect to the bundle
structure. Given local coordinates $(x,y)$ on $\pi^{-1}(\Omega)=\Omega\times\R^{k}$
and a smooth $n$-form $\omega\in E^{n}(\Omega\times\R^{k})$, $\omega:=\omega_{\alpha\beta}dx^{\alpha}\wedge dy^{\beta}$,
define auxiliary currents $E_{\alpha\beta}$ by $E_{\alpha\beta}(\omega):=\int\omega_{\alpha\beta}d\left\Vert T\right\Vert $,
where $\left\Vert T\right\Vert $ is the measure $\theta\mathcal{H}^{n}\rest Supp(T)$,
with $\mathcal{H}^{n}$ Hausdorff $n$-dimensional measure in $\Omega\times\R^{k}$
and $\theta$ the multiplicity of $T$ \cite[pp 45-46]{Morgan}. The
\emph{component currents} of $T$ are defined in terms of \emph{component
functions} $t_{\alpha\beta}:\Omega\times\R^{k}\rightarrow\R$ and
the auxiliary currents, by:\[
\left.T\right|_{\pi^{-1}(\Omega)}:=\left\{ T_{\alpha\beta}\right\} :=\left\{ t_{\alpha\beta}E_{\alpha\beta}\right\} .\]
The component functions $t_{\alpha\beta}:\pi^{-1}(\Omega)\rightarrow\R$
determine completely the current $T$, and the pairing between $T$
and an $n$-form $\omega\in E^{n}(B)\rest\Omega\times\R^{k}$ is given
by:\[
T(\omega):=\int_{\Omega\times\R^{k}}\sum_{\alpha\beta}t_{\alpha\beta}\omega_{\alpha\beta}d\left\Vert T\right\Vert .\]

\end{defn}

\begin{defn}
A bounded current $T$ in $B$ is \emph{a quasi-section} if, for each
coordinate neighborhood $\Omega\subset M$, 
\begin{enumerate}
\item $t_{n0}\geq0$ for $\left\Vert T\right\Vert $-almost all points $p\in Supp(T)$,
that is $<\overrightarrow{T}(q),\mathbf{e}(q)>\geq0$, $\left\Vert T\right\Vert $-almost
everywhere; where $\mathbf{e}(q):=\left.\frac{\partial}{\partial x^{1}}\wedge\cdots\wedge\frac{\partial}{\partial x^{n}}\right/\left\Vert \frac{\partial}{\partial x^{1}}\wedge\cdots\wedge\frac{\partial}{\partial x^{n}}\right\Vert $
is the (unique) horizontal $n$-plane at $q$ whose orientation is
preserved under $\pi_{*}$, and $\overrightarrow{T}$ is the unit
orienting $n$-vector field of $T$.
\item $\pi_{\#}(T)=1[M]$ as an $n$-dimensional current on $M$.
\item $\partial T=0$ (equivalently, for any $\Omega\subset M$, $\partial T\rest\pi^{-1}(\Omega)$
has support contained in $\partial\pi^{-1}(\Omega)$).
\end{enumerate}
There is an $M>0$ so that the fiber bundle $B$ is contained in the
disk bundle $E_{M}\subset E$ defined by $E_{M}:=\left\{ v\in E\left|\left\Vert v\right\Vert <M\right.\right\} $,
by compactness of $B$. Define the space $\widetilde{\Gamma}(E)$
to be the set of all countably rectifiable, integer multiplicity,
$n$-dimensional currents which are quasi-sections in $E$, with support
contained in $E_{M }$, called \emph{(bounded)} \emph{rectifiable
sections} of $E$. The space $\Gamma(E)$ of \emph{(strongly) rectifiable
sections of} $E$ is the smallest sequentially, weakly-closed space
containing the graphs of $C^{1}$ sections of $E$ which are supported
within $E_{M}$.
\end{defn}
A quasi-section which also is rectifiable is an element of $\widetilde{\Gamma}(E)$.
It would seem to be a strictly stronger condition to be in $\Gamma(E)$,
however, it is shown in \cite{Coventry} that, over a bounded domain
$\Omega$, $\Gamma(\Omega\times\R^{k})=\widetilde{\Gamma}(\Omega\times\R^{k})$.
The norm defined in \cite{Coventry} is finite in this case since
all currents have support contained in $E_{M}$. This extends to the
statement that $\Gamma(E)=\widetilde{\Gamma}(E)$ for a vector bundle
over a compact manifold $M$, since any such can be decomposed into
finitely many bounded domains where the bundle structure is trivial,
by a partition of unity argument.

The space $\widetilde{\Gamma}(B)$ of \emph{rectifiable sections}
of $B$ is the subset of $\widetilde{\Gamma}(E)$ of currents with
support in $B$, which is a weakly closed condition with respect to
weak convergence. This follows since, for any point $z$ outside of
$B$, there is a smooth form supported in a neighborhood of $z$ disjoint
from $B$. The space $\Gamma(B)$ of \emph{strongly} rectifiable sections
is the smallest sequentially, weakly-closed space containing the graphs
of $C^{1}$ sections of $B$. Since the fibers of $B$ are compact,
as is the base manifold $M$, minimal-mass elements will exist in
$\widetilde{\Gamma}(B)$ or $\Gamma(B)$, and mass-minimizing sequences
within any homology class will have convergent subsequences in $\widetilde{\Gamma}(B)$
or $\Gamma(B)$. This follows from lower semi-continuity with respect
to convergence of currents, convexity of the mass functional, and
the closure and compactness theorems for rectifiable currents. For
compact manifolds, as above, $\widetilde{\Gamma}(E)=\Gamma(E),$ but
it is \emph{not} the case that $\widetilde{\Gamma}(B)=\Gamma(B)$
in general. 

\begin{prop}
Let $\{ T_{j}\}\subset\widetilde{\Gamma}(B)$ (resp, $\Gamma(B)$)
be a sequence with equibounded mass. Then, there is a subsequence
which converges weakly to a current $T$ in $\widetilde{\Gamma}(B)$
(resp, $\Gamma(B)$) . 
\end{prop}
\begin{pf}
The Federer-Fleming compactness and closure theorems (see also \cite[Theorem 7.61, pp. 204-5]{Lin}
shows that a weak subsequence limit will exist and will be a countably-rectifiable,
integer-multiplicity current, with no interior boundaries. Since the
map $\pi:B\rightarrow M$ is proper, $\pi_{\#}$ will then commute
with weak limits, and so $\pi_{\#}(T)=1[M]$. Similarly, the conditions
$<\overrightarrow{T}(q),\mathbf{e}(q)>\geq0$ $\left\Vert T\right\Vert $-almost
everywhere and $Supp(T)\subset B$ are directly seen to be preserved
under weak limits, so the limit will be in $\widetilde{\Gamma}(B)$.
\hfill{}$\square$\end{pf}
\begin{defn}
Given a current $T$, the induced measures $\left\Vert T\right\Vert $
and $\left\Vert T_{\alpha\beta}\right\Vert $ are defined locally
by:\begin{eqnarray*}
\left\Vert T_{\alpha\beta}\right\Vert (A) & := & \sup\left(T_{\alpha\beta}(\omega)\right),\,\textrm{and}\\
\left\Vert T\right\Vert (A) & := & \sup\left(\sum_{\alpha\beta}T_{\alpha\beta}(\omega)\right),\end{eqnarray*}
where the supremum in either case is taken over all $n$-forms on
$B$, $\omega\in E_{0}^{n}(B)$, with $comass(\omega)\leq1$ \cite[ 4.1.7]{GMT}and
$Supp(\omega)\subset A$. 
\end{defn}

\section{Coordinatizability}

Let $T\in\widetilde{\Gamma}(B)$ have finite mass. Then, for each
$x\in M$, we say that $T$ is \emph{coordinatizable over} $x$ if
there is an $r>0$ so that $T\rest\pi^{-1}(B(x,r))$ (note that $\pi^{-1}(B(x,r))\cong B(x,r)\times F$)
has support contained within $B(x,r)\times U$, where $U\subset F$
is a contractible coordinate neighborhood of $F$, $U\cong\R^{j}$. 

\begin{prop}
The set of all points $x\in M$ where $T$ is coordinatizable over
$x$ is an open, dense subset of $M$.
\end{prop}
\begin{pf}
Openness follows from the definition, which involves open neighborhoods
in $M$. Note that the closed nested sets $Supp(T)\cap\pi^{-1}(\overline{B(x,r)})$,
as $r\rightarrow0$, have a nonempty intersection of $Supp(T)\cap\pi^{-1}(x)$.
So, given any neighborhood $U$ of $Supp(T)\cap\pi^{-1}(x)$ in $F$,
for some $r>0$ $\pi_{2}\left(Supp(T)\cap\pi^{-1}(B(x,r)\right)\subset U$,
where $\pi_{2}$ is the projection of $\pi^{-1}(B(x,r_{0}))\cong B(x,r_{0})\times F$
onto $F$. Certainly if $Supp(T)\cap\pi^{-1}(x)$ is finite, then,
since any finite set in $F$ is contained in a contractible coordinate
neighborhood in $F$, $T$ will be coordinatizable at $x$. So, any
point $x$ over which $T$ is not coordinatizable must have a preimage
under $\pi$ which is infinite, thus having infinite $0$-dimensional
Hausdorff measure. But, for\[
N:=\left\{ \left.x\in M\right|T\,\textrm{is}\, not\,\textrm{coordinatizable}\,\textrm{over}\,\textrm{x}\right\} ,\]
 then, if $N$ has positive Lebesgue measure on $M$, and if $\mathcal{F}$
is the volume (or mass) integrand,\begin{eqnarray*}
\mathcal{F}(T) & := & \int_{B}d\left\Vert T\right\Vert \\
 & \geq & \int_{\pi^{-1}(N)}d\left\Vert T\right\Vert \\
 & \geq & \int_{N}\#(\pi^{-1}(x))dx\\
 & = & \infty,\end{eqnarray*}
by the general area-coarea formula \cite{Morgan}.
\hfill{}$\square$\end{pf}
\begin{rem}
The Riemannian metric on $U\times V$ has the structure of a Riemannian
submersion $\pi:U\times V\rightarrow U$, that is, the projection
$\pi$ is an isometry on the orthogonal complement to the fibers,
and the projection onto the fiber, $\pi_{2}:U\times V\rightarrow V$
is an isometry restricted to each fiber. The fiber metric is not necessarily
Euclidean, and the orthogonal complements to the fibers will not necessarily
form an integrable distribution, but that will not affect the arguments
which follow.
\end{rem}

\section{Penalty Method}

Let $\mathcal{F}\,(=\mathcal{M})$ be the standard volume (area) functional,
applied to rectifiable sections. For an integer-multiplicity, countably-rectifiable
current $T=\tau(M,\theta,\overrightarrow{T}),$ where $M=Supp(T)$
and $\overrightarrow{T}$ is the unit orienting $n$-vector field,
as in \cite[p. 46]{Morgan}. Set, for each $\epsilon>0$, the modified
functional \[
\mathcal{F}_{\epsilon}(T):=\int_{T}f_{\epsilon}(\overrightarrow{T})d\left\Vert T\right\Vert ,\]
where $d\left\Vert T\right\Vert =\theta\mathcal{H}^{n}\rest Supp(T)$
and\[
f_{\epsilon}(\xi):=\left\Vert \xi\right\Vert +h_{\epsilon}(\xi):=\left\Vert \xi\right\Vert +\frac{1}{\epsilon}\left(\left|\xi_{n,0}\right|-\xi_{n,0}\right),\]
for $\xi\in\Lambda_{n}(T_{*}(B,z))\cong\Lambda_{n}(\R^{n+k})$ ($\left\Vert \xi\right\Vert $
is the usual norm of $\xi$ in $\Lambda_{n}(T_{*}(B,z))$ and $\xi_{n,0}:=<\xi,\mathbf{e>}$,
where $\mathbf{e}$ is the unique unit horizontal $n$-plane so that
$\pi_{*}(\mathbf{e)}=*dV_{M}$). 

Note also that, since the original integrand is positive, so is $f_{\epsilon}$,
at any point $\xi$. Moreover, $f_{\epsilon}$ satisfies the homogeneity
condition\[
f_{\epsilon}(t\xi)=tf_{\epsilon}(\xi)\]
 for $t>0$.

Set \[
\mathcal{H}_{0}(T):=\int_{T}h_{0}\left(\overrightarrow{T}\right)d\left\Vert T\right\Vert ,\]
where $h_{0}(\xi):=\left(\left|\xi_{n,0}\right|-\xi_{n,0}\right)$,
and set \[
\mathcal{H}_{\epsilon}(T):=\frac{1}{\epsilon}\mathcal{H}_{0}(T).\]

On the parts of $T$ which project to a negatively-oriented current
(locally) on the base, the functional $\mathcal{H}_{0}()$ has value
equal to twice the Lebesgue measure of the projected image, considered
as measurable subsets of the base. 

Clearly $f_{\epsilon}$ satisfies the bounds \[
\left\Vert \xi\right\Vert \leq f_{\epsilon}(\xi)\leq\left(1+\frac{2}{\epsilon}\right)\left\Vert \xi\right\Vert .\]
 In addition, the functional satisfies the $\lambda$-ellipticity
condition with $\lambda=1$ \begin{equation}
\left[\mathcal{M}(X)-\mathcal{M}(mD)\right]\leq\mathcal{F}_{\epsilon}(X)-\mathcal{F}_{\epsilon}(mD)\label{eq:ellipticity}\end{equation}
 where $mD$ is a flat disk with multiplicity $m$ and $X$ is a rectifiable
current with the same boundary as $mD$. This inequality is clear
if $\pi_{\#}(mD)$ is positively-oriented, since in that case $\mathcal{M}(mD)=\mathcal{F}_{\epsilon}(mD)$,
and (in all cases) $\mathcal{M}(X)\leq\mathcal{F}_{\epsilon}(X)$.
If $\pi_{\#}(mD)$ is negatively-oriented, though, then $\pi_{\#}(X)=\pi_{\#}(mD)$
since they have the same boundary and are integer-multiplicity countably-rectifiable
$n$-currents on $\R^{n}$, by the constancy theorem. However, in
this case $\mathcal{H}_{\epsilon}(X)\geq\mathcal{H}_{\epsilon}(mD)$,
and the result follows.

\subsection{Minimization problem}

We now consider the mass-minimization problem for rectifiable sections
$T\in\widetilde{\Gamma}(B)$ within a given integral homology class
$[T]\in H_{n}(B,\mathbb{Z})$ which includes graphs, that is, for
which there is a smooth section $S_{0}\in\Gamma(B)$ with $[S_{0}]=[T]$.
Set $A:=\left\Vert S_{0}\right\Vert $. Set \[
\mathbf{R}[T]:=\left\{ \left.S\in[T]\right|S\textrm{ is a countably rectifiable, integer-multiplicity $n$-current in }B\right\} .\]

For any $\epsilon>0$, since the tangent planes at each point of $S_{0}$
projects to an $n$-plane of positive orientation, $h_{\epsilon}(\overrightarrow{S_{0}})=\frac{1}{\epsilon}\left(\left|\xi_{n,0}\right|-\xi_{n,0}\right)=0$,
and so $\mathcal{F}_{\epsilon}(S_{0})=\left\Vert S_{0}\right\Vert :=A$,
which shows that $\left\{ S\in\mathbf{R}[T]\left|\mathcal{F}_{\epsilon}(S)\leq A\right.\right\} \neq\emptyset$.
Thus, if $B_{0}:=\left\{ S\in\mathbf{R}[T]\left|\left\Vert S\right\Vert \leq2A\right.\right\} $,
\[
Lev_{A}\mathcal{F}_{\epsilon}:=\left\{ S\in\mathbf{R}[T]\left|\mathcal{F}_{\epsilon}(S)\leq A\right.\right\} \subset B_{0},\]
since for any current $\mathcal{F}_{\epsilon}(S)\geq\left\Vert S\right\Vert $.
Also, by the Federer-Fleming closure theorem, $B_{0}$ is compact
with respect to the usual convergence of currents. Since the functional
$\mathcal{F}_{\epsilon}$ is elliptic (eq (\ref{eq:ellipticity})),
it will be lower semi-continuous with respect to weak convergence
of rectifiable currents \cite[ 5.1.5]{GMT}. Thus each $Lev_{A}\mathcal{F}_{\epsilon}$
is compact in this topology, and so by \cite{Struwe}, for each such
$\epsilon$, an $\mathcal{F}_{\epsilon}$-energy-minimizing rectifiable
current $T_{\epsilon}\in[T]$ exists, and $\mathcal{F}_{\epsilon}(T_{\epsilon})<\left\Vert S_{0}\right\Vert =A$. 

Set \begin{eqnarray*}
min(\mathcal{F}_{\epsilon}) & := & min\left\{ \mathcal{F}_{\epsilon}(T)\left|T\in\mathbf{R}[T]\right.\right\} \\
Argmin(\mathcal{F}_{\epsilon}) & := & \left\{ T\in\mathbf{R}[T]\left|\mathcal{F}_{\epsilon}(T)=min(\mathcal{F}_{\epsilon})\right.\right\} ,\\
min(\mathcal{F}) & := & min\left\{ \left.\mathcal{F}(T)\right|T\in[T]\cap\Gamma(B)\right\} ,\,\,\textrm{and}\\
Argmin(\mathcal{F}) & := & \left\{ T\in[T]\cap\Gamma(B)\left|\mathcal{F}(T)=min(\mathcal{F})\right.\right\} .\end{eqnarray*}
 Similarly to \cite{Roubicek}, we have

\begin{prop}
\textbf{\emph{{[}Convergence of the penalty problems{]}}} \label{thm:penalty argument}\begin{eqnarray*}
\lim_{\epsilon\downarrow0}min\left(\mathcal{F}_{\epsilon}\right) & = & min\left(\mathcal{F}\right),\\
\limsup_{\epsilon\downarrow0}Argmin(\mathcal{F}_{\epsilon}) & \subset & Argmin(\mathcal{F}).\end{eqnarray*}

\end{prop}
\begin{rem}
That is, the minimal values of the penalty functionals on that homology
class converge to the minimum of the mass of all homologous rectifiable
sections, and the limsup of the set of minimizing currents \cite{Roubicek}
of the penalty problems is contained in the set of mass-minimizing
rectifiable sections. This does not imply that each mass-minimizing
rectifiable section is the limit of a sequence of minimizers of the
penalty problems, but that one such mass-minimizing rectifiable section
is such a limit.
\end{rem}
\begin{pf}
Since the set of countably-rectifiable integer-multiplicity currents
in $[T]$ (the domain of $\mathcal{F}_{\epsilon}$) contains the rectifiable
sections, and $\mathcal{F}_{\epsilon}(S)=\mathcal{F}(S)=\left\Vert S\right\Vert $
for any rectifiable section $S$, we have immediately that $min\left(\mathcal{F}_{\epsilon}\right)\leq min\left(\mathcal{F}\right)$.
Moreover, $min\left(\mathcal{F}_{\epsilon_{1}}\right)\leq min\left(\mathcal{F}_{\epsilon_{2}}\right)$
if $\epsilon_{1}>\epsilon_{2}$, since for $T_{\epsilon_{i}}$ minimizers
of $\mathcal{F}_{\epsilon_{i}}$,\[
\mathcal{F}_{\epsilon_{1}}(T_{\epsilon_{1}})\leq\mathcal{F}_{\epsilon_{1}}(T_{\epsilon_{2}})\leq\mathcal{F}_{\epsilon_{2}}(T_{\epsilon_{2}}),\]
so $\lim_{\epsilon\downarrow0}min(\mathcal{F}_{\epsilon})$ exists. 

Take some $T_{\epsilon}\in\mathbf{R}[T]$ which minimizes $\mathcal{F}_{\epsilon}$within
$\mathbf{R}[T]$. Then \begin{eqnarray*}
\left\Vert T_{\epsilon}\right\Vert  & = & \mathcal{F}_{\epsilon}(T_{\epsilon})-H_{\epsilon}(T_{\epsilon})\\
 & \leq & \mathcal{F}_{\epsilon}(T_{\epsilon})\\
 & = & min(\mathcal{F}_{\epsilon})\\
 & \leq & min(\mathcal{F}).\end{eqnarray*}
This shows that $T_{\epsilon}\in B_{0}$ above, which, in the topology
of weak convergence of countably-rectifiable, integer-multiplicity
currents, is compact. So, by the Federer-Fleming compactness and closure
theorems \cite[4.2.16, 4.2.17]{GMT}, some subsequence of $\left\{ T_{\epsilon}\right\} $
converges as $\epsilon\downarrow0$ to some $S\in B_{0}$. 

Since the penalty component satisfies\[
\mathcal{H}_{0}(T_{\epsilon})=\epsilon\mathcal{H}_{\epsilon}(T_{\epsilon})=\epsilon\left(min(\mathcal{F}_{\epsilon})-\left\Vert T_{\epsilon}\right\Vert \right),\]
and the penalty functional $\mathcal{F}_{\epsilon}$ is lower semi-continuous
with respect to weak convergence of currents, we have \begin{eqnarray*}
\mathcal{H}_{0}(S) & \leq & \liminf_{\epsilon\downarrow0}(\mathcal{H}_{0}(T_{\epsilon}))\\
 & = & \liminf_{\epsilon\downarrow0}\epsilon\left(min(\mathcal{F}_{\epsilon})-\left\Vert T_{\epsilon}\right\Vert \right)\\
 & \leq & \liminf_{\epsilon\downarrow0}\epsilon\left(min(\mathcal{F})\right)\\
 & = & 0.\end{eqnarray*}
So. immediately we have that $S\in\widetilde{\Gamma}(B)$, so that
$\mathcal{F}(S)\geq min(\mathcal{F})$. Applying the same limit to
the previous equation, \begin{eqnarray*}
\mathcal{F}(S) & = & \left\Vert S\right\Vert \\
 & \leq & \liminf_{\epsilon\downarrow0}\left\Vert T_{\epsilon}\right\Vert \\
 & = & \liminf_{\epsilon\downarrow0}\left(\mathcal{F}_{\epsilon}(T_{\epsilon})-H_{\epsilon}(T_{\epsilon})\right)\\
 & \leq & \liminf_{\epsilon\downarrow0}\mathcal{F}_{\epsilon}(T_{\epsilon})\\
 & = & \liminf_{\epsilon\downarrow0}min(\mathcal{F}_{\epsilon})\\
 & \leq & min(\mathcal{F}),\end{eqnarray*}
which implies that all inequalities must be equalities, and $S$ is
a mass-minimizing rectifiable section in $[T]\cap\widetilde{\Gamma}(B)$.
In addition, we get that \[
\lim_{\epsilon\downarrow0}min(\mathcal{F}_{\epsilon})=min(\mathcal{F})\]
and, any limit current of a subsequence of minimizers $\left\{ T_{\epsilon}\right\} $(for
a sequence of $\epsilon$'s going to 0) will be a minimizer $T_{0}$
of \emph{$\mathcal{F}$ on $[T]\cap\widetilde{\Gamma}(B)$.}
\hfill{}$\square$\end{pf}
The set of points $B_{\epsilon}\subset\Omega$ where $T_{\epsilon}$
is not a section, \[
B_{\epsilon}:=\pi\left(\left\{ \left.p\in Supp(T_{\epsilon})\right|h_{\epsilon}(\overrightarrow{T}_{p})>0\right\} \right),\]
satisfies, where $\mathbf{e}$ is the horizontal $n$-plane, \begin{eqnarray}
\mathcal{H}_{0}(T_{\epsilon}) & = & \epsilon\mathcal{H}_{\epsilon}(T_{\epsilon})\nonumber \\
 & = & \int_{\Omega\times F}h_{0}(T_{\epsilon})d\left\Vert T_{\epsilon}\right\Vert \nonumber \\
 & = & \int_{\Omega\times F}(\left|<\overrightarrow{T_{\epsilon}},\mathbf{e}>\right|-<\overrightarrow{T_{\epsilon}},\mathbf{e}>)d\left\Vert T_{\epsilon}\right\Vert \nonumber \\
 & = & -2\int_{B_{\epsilon}\times F}<\overrightarrow{T_{\epsilon}},\mathbf{e}>d\left\Vert T_{\epsilon}\right\Vert \nonumber \\
 & = & 2\int_{B_{\epsilon}\times F}\left|<\overrightarrow{T_{\epsilon}},\mathbf{e}>\right|d\left\Vert T_{\epsilon}\right\Vert \nonumber \\
 & = & 2T_{\epsilon}\left(\pi^{*}(d\left.V_{\Omega}\right|_{B_{\epsilon}})\right)\nonumber \\
 & = & 2\pi_{\#}(T_{\epsilon})(d\left.V_{\Omega}\right|_{B_{\epsilon}})\nonumber \\
 & \geq & 2\left\Vert B_{\epsilon}\right\Vert ,\label{eq:B-epsilon}\end{eqnarray}
where $dV_{\Omega}$ is the volume element of the base, since $\left.\pi_{\#}(T_{\epsilon})\right|_{B_{\epsilon}}$
is a (positive integer) multiple of the fundamental class of the base,
restricted to $B_{\epsilon}$. From the previous result,\[
\lim_{\epsilon\downarrow0}\mathcal{H}_{\epsilon}(T_{\epsilon})=0,\]
thus $\left\Vert B_{\epsilon}\right\Vert $ approaches 0 more rapidly
than $\epsilon$ itself. Similarly to \cite{Struwe}, we have the
following:

\begin{lem}
\label{lem:bad-set-size} $\left\Vert B_{\epsilon}\right\Vert \leq\frac{\epsilon}{\left|\ln(\epsilon)\right|}A$,
where $A$ depends only on dimension and the homology class $[T]\in H_{n}(B,\mathbb{Z})$.
\end{lem}
\begin{pf}
If $v_{\epsilon}:=\mathcal{F}_{\epsilon}(T_{\epsilon})$, for $0<\epsilon_{1}<1$,
then since $v_{\epsilon}$ is a monotone-decreasing function of $\epsilon$,
it is differentiable almost-everywhere, and \begin{eqnarray*}
\left|v_{\epsilon}'\right| & = & \left|\lim_{h\rightarrow0}\frac{v_{\epsilon}-v_{\epsilon-h}}{h}\right|\\
 & \geq & \left|\lim_{h\rightarrow0}\frac{\mathcal{F}_{\epsilon}(T_{\epsilon})-\mathcal{F}_{\epsilon-h}(T_{\epsilon})}{h}\right|\\
 & = & \left|\lim_{h\rightarrow0}\frac{\left(\frac{1}{\epsilon}-\frac{1}{\epsilon-h}\right)}{h}\right|\mathcal{H}_{0}(T_{\epsilon})\\
 & = & \frac{1}{\epsilon^{2}}\mathcal{H}_{0}(T_{\epsilon}).\end{eqnarray*}
 In addition, for any fixed rectifiable section $S$ in the homology
class $[T]$, for all $\epsilon>0$ $\nu_{\epsilon}\leq\mathcal{F}(S)$,
so that $\nu_{\epsilon}$ is bounded. 

Now, as in \cite[p. 70, Theorem 7.3]{Struwe}, \begin{eqnarray*}
C & \ge & v_{\epsilon_{1}}-v_{1}\\
 & \geq & \int_{\epsilon_{1}}^{1}\left|v_{\epsilon}'\right|d\epsilon\\
 & \geq & \int_{\epsilon_{1}}^{1}\textrm{ess $\inf_{\epsilon_{1}<\epsilon<1}$}\left(\epsilon\left|v_{\epsilon}'\right|\right)\frac{1}{\epsilon}d\epsilon\\
 & = & \textrm{ess $\inf_{\epsilon_{1}<\epsilon<1}$}\epsilon\left|v_{\epsilon}'\right|\cdot\left(-\ln(\epsilon_{1})\right)\\
 & \geq & \textrm{ess $\inf_{\epsilon_{1}<\epsilon<1}$}\frac{1}{\epsilon}\mathcal{H}_{0}(T_{\epsilon})\cdot\left|\ln(\epsilon_{1})\right|\\
 & \geq & \textrm{ess $\inf_{\epsilon_{1}<\epsilon<1}$}\frac{1}{\epsilon}\left\Vert B_{\epsilon}\right\Vert \cdot\left|\ln(\epsilon_{1})\right|.,\end{eqnarray*}
applying (\ref{eq:B-epsilon}). Since $v_{\epsilon_{1}}-v_{1}$is
bounded (and nonnegative), there is a constant $C$ so that \[
\left\Vert B_{\epsilon}\right\Vert \leq\frac{\epsilon}{\left|\log(\epsilon)\right|}C,\]
where $C$ depends only on the homology class $[T]$ of sections being
considered. 
\hfill{}$\square$\end{pf}

\section{Existence of tangent cones}

Let $T$ be a mass-minimizing rectifiable section, and presume that
$T$ is the limit of a sequence $T_{\epsilon_{i}}$ of minimizers
of the penalty energy $\mathcal{F}_{\epsilon_{i}}$. (At least one
minimizer of the mass functional among rectifiable sections is of
this form), by Proposition (\ref{thm:penalty argument}).

\begin{prop}
\label{pro:monotonicity}For any point $p\in Supp(T)$, the mass-density
$\Theta(p,T)$ is at least 1. Moreover, there is a (possibly non-unique)
tangent cone at $p$ of $T$.
\end{prop}
\begin{rem}
The proof will depend on a monotonicity of mass ratio result. Once
that is established, the result will follow similarly to the case
for area-minimizing rectifiable currents.
\end{rem}
\begin{lem}
\textbf{\emph{{[}Monotonicity of mass ratio{]}.}} For any $p\in Supp(T)$,
the ratio\[
\frac{\mathcal{F}\left(T\rest B(p,r)\right)}{r^{n}}\]
is a monotone increasing function of $r$.
\end{lem}
\begin{pf}
(of the Lemma). Consider, for a sequence $\epsilon=\epsilon_{i}$
converging to $0$, the penalty energy function\[
f_{\epsilon}(r):=\mathcal{F}_{\epsilon}\left(T_{\epsilon}\rest B(p_{\epsilon},r)\right),\]
 where $p_{\epsilon}\in Supp(T_{\epsilon})$. We show that the penalty
function satisfies the monotonicity differential inequality $(f_{\epsilon}(r)/r^{n})'\geq0$,
as in \cite{Morgan}. 

Choose a radius $r$ for which the boundary $\partial\left(T_{\epsilon}\rest B(p_{\epsilon},r)\right)$
is rectifiable (true for almost-all $r$ by slicing). For such an
$r$, note that $\partial(T_{\epsilon}\rest B(p_{\epsilon},r))$ is
the boundary of the restriction of $T_{\epsilon}$ to the ball. Let
$C[\partial(T_{\epsilon}\rest B(p_{\epsilon},r))]$ be the cone over
$\partial(T_{\epsilon}\rest B(p_{\epsilon},r))$ with cone point $p_{\epsilon}$,
oriented so that $C[\partial(T_{\epsilon}\rest B(p_{\epsilon},r))]+T_{\epsilon}\rest(B\backslash B(p_{\epsilon},r))$
is a cycle. Define a boundary penalty-energy $\partial\mathcal{F}_{\epsilon}$
by restriction, that is: \[
\partial\mathcal{F}_{\epsilon}(\partial(T_{\epsilon}\rest B(p_{\epsilon},r))):=\int_{B}\left\Vert \overrightarrow{T_{\epsilon}}\right\Vert +\frac{1}{\epsilon}\left(\left|<\overrightarrow{T_{\epsilon}},\mathbf{e}>\right|-<\overrightarrow{T_{\epsilon}},\mathbf{e}>\right)d\left\Vert \partial(T_{\epsilon}\rest B(p_{\epsilon},r))\right\Vert .\]
Since $C_{r}:=C[\partial((T_{\epsilon}\rest G_{\epsilon})\rest B(p_{\epsilon},r))]$
is a cone, \begin{eqnarray*}
\mathcal{F}_{\epsilon}(C_{r}) & \leq & \frac{n}{r}\partial\mathcal{F}_{\epsilon}(\partial C_{r})\\
 & = & \frac{n}{r}\partial\mathcal{F}_{\epsilon}(\partial(T_{\epsilon}\rest B(p_{\epsilon},r))).\end{eqnarray*}
Now, set \[
f_{\epsilon}(r):=\mathcal{F}_{\epsilon}(T_{\epsilon}\rest B(p_{\epsilon},r)).\]
We claim that slicing by $u(x)=\left\Vert x-p_{\epsilon}\right\Vert $
yields that, for almost-every $r$ (as above)\[
\partial\mathcal{F}_{\epsilon}(\partial(T_{\epsilon}\rest B(p_{\epsilon},r)))\leq f_{\epsilon}'(r).\]
To show this, let $T$ be a rectifiable current, and $u$ Lipschitz.
The slice \[
<T,u,r+>:=\partial T\rest\left\{ x\left|u(x)>r\right.\right\} -\partial(T\rest\left\{ x\left|u(x)>r\right.\right\} )\]
 satisfies, for $\partial\mathcal{H}_{0}(<T,u,r+>):=\int_{B}\left(\left|<\overrightarrow{T},\mathbf{e}>\right|-<\overrightarrow{T},\mathbf{e}>\right)d\left\Vert <T,u,r+>\right\Vert $,
the following:\begin{eqnarray*}
\partial\mathcal{H}_{0}(<T,u,r+>) & \leq & Lip(u)\liminf_{h\downarrow0}\mathcal{H}_{0}(T)\rest\{ r<u(x)<r+h\}/h\\
 & = & Lip(u)\frac{\partial}{\partial r}\mathcal{H}_{0}(T)\rest\left\{ x\left|u(x)\leq r\right.\right\} ,\end{eqnarray*}
where we have abused notation and denoted the Dini derivative in the
previous line by $\partial/\partial r$. This follows by considering,
for a small, positive $h$, a smooth approximation $f$ of the characteristic
function of $\left\{ x\left|u(x)>r\right.\right\} $ with 
\[
f(x)=\left\{ \begin{array}{cc}
0, & \textrm{if }u(x)\leq r\\
1, & \textrm{if }u(x)\geq r+h\end{array}\right.\]
and $Lip(f)\leq Lip(u)/h$. Then (cf. \cite[4.11, p. 56]{Morgan})\begin{eqnarray*}
\partial\mathcal{H}_{0}(<T,u,r+>) & \approx & \partial\mathcal{H}_{0}((\partial T)\rest f-\partial(T\rest f))\\
 & = & \partial\mathcal{H}_{0}(T\rest df)\\
 & \leq & Lip(f)\mathcal{H}_{0}(T)\rest\{ r<u(x)<r+h\}\\
 & \lesssim & Lip(u)\mathcal{H}_{0}(T)\rest\{ r<u(x)<r+h\}/h\\
 & = & Lip(u)\frac{\partial}{\partial u}\mathcal{H}_{0}(T)\rest\left\{ x\left|u(x)\leq r\right.\right\} .\end{eqnarray*}
In the present case, with $u(x):=\left\Vert x-p_{\epsilon}\right\Vert $,
$<T_{\epsilon},u,r+>=\partial(T_{\epsilon}\rest B(p_{\epsilon},r))$,
$\partial\mathcal{F}_{\epsilon}(\partial(T_{\epsilon}\rest B(p_{\epsilon},r)))\leq f_{\epsilon}'(r)$
as claimed for almost-every $r$, since for the standard mass functional
this result is standard, and $\mathcal{F}_{\epsilon}=\mathcal{M}+\frac{1}{\epsilon}\mathcal{H}_{0}$. 

Combining these two relationships together and using minimality of
$T_{\epsilon}$,\[
f_{\epsilon}(r):=\mathcal{F}_{\epsilon}(T_{\epsilon}\rest B(p_{\epsilon},r))\leq\mathcal{F}_{\epsilon}(C[\partial(T_{\epsilon}\rest B(p_{\epsilon},r))])\leq\frac{n}{r}\partial\mathcal{F}_{\epsilon}(\partial(T_{\epsilon}\rest B(p_{\epsilon},r)))\leq\frac{n}{r}\frac{df_{\epsilon}(r)}{dr},\]
for almost-every $r$, hence the absolutely continuous part of $f_{\epsilon}(r)/r^{n}$
is increasing. Since any singular part is due to increases in $f_{\epsilon}(r)$,
$f_{\epsilon}(r)/r^{n}$ is increasing as claimed.

Let $p_{\epsilon}\rightarrow p$ be a sequence of points on the support
of the penalty minimizers converging to $p\in Supp(T)$. Set $f(r):=\mathcal{F}\left(T\rest B(p,r)\right)$.
Since $f_{\epsilon}(r)/r^{n}$ is monotone increasing as a function
of $r$ for each fixed $\epsilon>0$, so will be $f(r)/r^{n}$.
\hfill{}$\square$\end{pf}

Arguing precisely as in \cite[ 5.4.3]{GMT}, (see also \cite[pp. 90-95]{Morgan}),
Proposition (\ref{pro:monotonicity}) follows.

\section{Domain of the penalty-minimizers}

Let $\Omega=B(x_{0},R)$ be a ball. It follows from the structure
theorem for rectifiable currents \cite{GMT} that, except over the
bad set $B_{\epsilon}$, which is a set of mass less than $\epsilon R^{n}$,
the penalty-minimizer $T_{\epsilon}$ will be the graph of a vector-valued
BV function $u_{\epsilon}$.

The points $x\in\Omega\backslash B_{\epsilon}$ so that for all $p=\pi^{-1}(x)\cap Supp(T_{\epsilon})$,
$\Theta(p)=1$, has measure approaching that of $\Omega$ as $\epsilon\downarrow0$
because of our bounds on $B_{\epsilon}$. Since $T_{\epsilon}\rest\pi^{-1}(\Omega\backslash B_{\epsilon})$
is a rectifiable section, the structure theorem for rectifiable currents
implies that for $\Omega$-a.e. points $x$ of $\Omega\backslash B_{\epsilon}$,
there is one point in $\pi^{-1}(x)\cap Supp(T_{\epsilon})$. 

Define $u_{\epsilon}$ as a vector-valued BV-function over $\Omega\backslash B_{\epsilon}$
whose carrier is $Supp(T_{\epsilon})$ \cite[Section IV]{Bombieri},
defined coordinatewise by integration, first defining $S_{j}$ as
$n$-dimensional currents in $U$ by $S_{j}(\phi):=T(y_{j}\pi^{*}(\phi))$
for $\phi\in E^{n}(U)$ and $y_{j}$ the $j^{\textrm{th}}$ coordinate
of the fiber ($U$ must be a coordinatizable neighborhood). Then,
the components of $u_{\epsilon}$ can be defined by $S_{j}(\phi)=\int(u_{\epsilon})_{j}(x)\phi$,
which define the components as $\textrm{BV}_{\textrm{loc}}$-functions
on $U$. 

It is not clear (compare \cite[p. 106]{Bombieri}) that this BV map
will be a Lipschitz graph a.e. in general. For example, if $T$ is
the simple staircase current $T_{\alpha}=[[(t,\alpha\left\lfloor t\right\rfloor )]]+[[(\left\lfloor t\right\rfloor ,\alpha(t-1))]],\, t\in[0,n]$,
$T\in\Gamma([0,n]\times\R)$, then $T$ will be a polyhedral chain,
and so the image of a Lipschitz map. However, the set $A$ on which
$T\rest\pi^{-1}(A)$ will have a single point in each preimage is
the base interval minus finitely many points (excluding the points
that are the projections of the risers of the stairs), and $Supp(T)\cap\pi^{-1}(A)$
cannot be a Lipschitz graph on all of $A$. By controlling the height
$\alpha$ of the risers the total cylindrical excess $E$ of this
example can be as small as needed as well.

However, it is the case that there will be, for any positive number
$\delta>0$, a Lipschitz map $g$ so that $g=u_{\epsilon}$ except
on a set of measure less than $\delta$, by Theorem 2 page 252 of
\cite{Evans and Gariepy}. In fact, $g$ can be taken to be $C^{1}$
by Corollary 1, p. 254, of the same reference. The Lipschitz constant
of the map $g$ will clearly depend upon $\delta$, as is illustrated
by the example above.

Now, it is not necessarily true that the graph of $g$ will agree
with $Supp(T_{\epsilon})$ on the set where $g$ agrees with $u_{\epsilon}$,
since that graph does not necessarily agree with $Supp(T_{\epsilon})$
itself.

\begin{prop}
\label{pro:g}For any $\epsilon>0$, there is a set $D_{\epsilon}\supseteq B_{\epsilon}$
of measure less than $2\left\Vert B_{\epsilon}\right\Vert $ and a
$C^{1}$ map $g_{\epsilon}:U\backslash D_{\epsilon}\rightarrow F$
so that, as rectifiable currents, \[
graph(g_{\epsilon})\rest\pi^{-1}(U\backslash D_{\epsilon})=T_{\epsilon}\rest\pi^{-1}(U\backslash D_{\epsilon}).\]

\end{prop}
\begin{pf}
For $\delta>0$ sufficiently small, Choose $g_{\epsilon}$ by \cite[Corollary 1, p. 254]{Evans and Gariepy}
to agree with $u_{\epsilon}$ on $U$ except for a set of measure
$\delta$, and to be $C^{1}$ and Lipschitz there. Take $D_{\epsilon}$
to be the union of this set with $B_{\epsilon}$, which if $\delta$
is chosen small enough will have measure bounded by $2\left\Vert B_{\epsilon}\right\Vert $.
It suffices to show that these currents agree except over a set of
measure 0 in the domain, outside of $D_{\epsilon}$. However, if they
disagree on a set $A$, within $U\backslash D_{\epsilon}$, of positive
measure, then for some $i$, $g_{j}(x)=(u_{\epsilon})_{j}(x)$ is
different from $y_{j}(Supp(T_{\epsilon})\cap\pi^{-1}(x))$ on $A$.
However, for any form $\phi$ on the base, over any subset $V\subseteq U\backslash D_{\epsilon}$,
since $\pi_{\#}(T)=1\cdot[U]$, \begin{eqnarray*}
 &  & \int_{V}y_{j}(Supp(T_{\epsilon})\cap\pi^{-1}(x))\,\phi\\
 & = & \int_{Supp(T_{\epsilon})\cap\pi^{-1}(V)}y_{j}(Supp(T_{\epsilon})\cap\pi^{-1}(x))<\pi^{*}(\phi),\overrightarrow{T}_{\epsilon}>d\mathcal{H}^{n}\\
 & = & \int_{Supp(T_{\epsilon})\cap\pi^{-1}(V)}<y_{j}\pi^{*}(\phi),T_{\epsilon}>d\mathcal{H}^{n}\\
 & = & \left(T_{\epsilon}\rest\pi^{-1}(V)\right)(y_{j}\pi^{*}(\phi))\\
 & = & S_{j}(\phi)\\
 & = & \int_{V}(u_{\epsilon})_{j}(x)\phi.\end{eqnarray*}
Since this equality must hold for all $\phi$ and $V\subset U\backslash D_{\epsilon}$
as above, the two functions must agree on a set of full measure.
\hfill{}$\square$\end{pf}
\begin{note}
The mass $\left\Vert D_{\epsilon}\right\Vert $ will satisfy \[
\lim_{\epsilon\rightarrow0}\frac{1}{\epsilon}\left\Vert D_{\epsilon}\right\Vert =0\]
 by the construction of both $B_{\epsilon}$ and the extension $D_{\epsilon}$
as defined in the proof of the previous result. Similarly, Lemma {[}\ref{lem:bad-set-size}{]}
will imply that \[
\left\Vert D_{\epsilon}\right\Vert \leq\frac{2\epsilon}{\left|\log(\epsilon)\right|}A,\]
with $A$ depending only on dimension.
\end{note}

\section{Homotopies and deformations}

Let $T^{t}$ be a one-parameter family of countably-rectifiable integer-multiplicity
currents with $T^{0}=T_{\epsilon}$, smooth in $t$. The derivative
$h:=\left.\frac{d}{dt}\right|_{0}T^{t}$ at $t=0$ is a current, but
in general will not be a rectifiable current. The support of $h$
will be $T_{\epsilon}$, but $h$ will be represented by integration
as \[
h(\phi):=\int_{E}<\overrightarrow{h},\phi>d\left\Vert T_{\epsilon}\right\Vert ,\]
where \[
\overrightarrow{h}d\left\Vert T_{\epsilon}\right\Vert =\left.\frac{d}{dt}\right|_{0}\overrightarrow{T^{t}}d\left\Vert T^{t}\right\Vert .\]
If $T_{\epsilon}$ is a smooth graph, $T_{\epsilon}=graph(g_{\epsilon})$,
then $T^{t}$ will be also, for $t$ sufficiently small, $T^{t}=graph(g_{\epsilon}+tk+\mathcal{O}(t^{2}))$,
by the implicit function theorem, and \begin{eqnarray*}
\overrightarrow{h}d\left\Vert T_{\epsilon}\right\Vert  & = & \left.\frac{d}{dt}\right|_{0}\overrightarrow{T^{t}}d\left\Vert T^{t}\right\Vert \\
 & = & \left.\frac{d}{dt}\right|_{0}(e_{1}+\nabla_{1}g_{\epsilon}+t\nabla_{1}k+t^{2}*)\wedge\cdots\wedge(e_{n}+\nabla_{n}g_{\epsilon}+t\nabla_{n}k+t^{2}*)d\mathcal{L}^{n}\\
 & = & \left(\nabla_{1}k\wedge(e_{2}+\nabla_{2}g_{\epsilon})\wedge\cdots\wedge(e_{n}+\nabla_{n}g_{\epsilon})+\cdots+(e_{1}+\nabla_{1}g_{\epsilon})\wedge\cdots\wedge\nabla_{n}k\right)d\mathcal{L}^{n}.\end{eqnarray*}

\begin{rem}
\label{rem:first-order-in-dh}Note that this derivative is first-order
with respect to the derivative $Dk$. The derivative will be first-order
with respect to $Dk$ for places where $T_{\epsilon}$ is not a graph,
since, being rectifiable, $Dk$ is a sum of terms of that sort. 
\end{rem}
Equivalently, we can consider maps $H_{t}:[0,1]\times U\times\R^{j}\rightarrow U\times\R^{j}$,
ambient homotopies of the region into itself, and the push-forward
$(H_{t})_{\#}(T)=T^{t}$. Of particular interest will be in families
which are \emph{vertical} in the sense that $H_{t}(x,y)=(x,y+\eta(t,x))$
for some $\eta:[0,1]\times U\rightarrow\R^{j}$. These are, of course,
in the graph case equivalent to families $T^{t}=graph(g_{\epsilon}+\eta(t,x))$.

\subsection{\label{sub:Euler-Lagrange-equations-for-T-epsilon}Euler-Lagrange
equations for $T_{\epsilon}$}

Restrict the deformations $T^{t}$ to be, for each $\epsilon>0$,
deformations in the vertical directions only. For a rectifiable section,
such a deformation will remain a section. If the domain $U=B(x_{0},R),$
is a coordinatizable neighborhood, so that the fiber can be considered
to be a compact subset of $\R^{j}$, and if coordinates are chosen
so that $(x_{0},\overline{y})$ is $(0,0)$ (for a particular value
of $\overline{y}$ to be determined), then, following \cite{Bombieri},
a deformation given by $T^{t}=(H_{t,R})_{\# }(T_{\epsilon})$, where
\begin{equation}
H_{t,R}(x,y)=(x,y+t\eta(x/R)),\label{eq:homotopy}\end{equation}
so that, over $B(x_{0},R)\backslash D_{\epsilon}$, $T_{\epsilon}\rest C(x_{0},R)=graph(g_{t})$,
where $g_{t}(x)=g_{\epsilon}(x)+t\eta(x/R)$, and where $\eta:B(0,1)\rightarrow\R^{k}$
is a smooth test function with support within the open ball and with
$\left\Vert \nabla\eta\right\Vert \leq1$ pointwise. Set $H_{t}:=H_{t,1}$.

Over a set of full measure in $Supp(T_{\epsilon})$ the tangent cone
at $(x,g_{\epsilon}(x))\in Supp(T_{\epsilon})$ is an $n$-plane and
is defined as usual from the graph of $g_{\epsilon}$. Since the area
functional, as a functional over the base, is then \[
\int_{\Omega\backslash D_{\epsilon}}\sqrt{1+\left\Vert \nabla g_{\epsilon}\right\Vert ^{2}+\left\Vert \nabla g_{\epsilon}\wedge\nabla g_{\epsilon}\right\Vert ^{2}+\cdots+\left\Vert \nabla g_{\epsilon}\begin{array}{c}
n\\
\overbrace{\wedge\cdots\wedge}\end{array}\nabla g_{\epsilon}\right\Vert ^{2}}d\mathcal{L}^{n},\]
then the Euler-Lagrange equations, obtained from calculus of variations
methods (using a vertical variation $g_{t}(x)=g_{\epsilon}(x)+t\eta(x/R)$),
is\begin{eqnarray*}
 &  & \left.\frac{d}{dt}\right|_{0}\int_{\Omega\backslash D_{\epsilon}}\sqrt{1+\left\Vert \nabla g_{t}\right\Vert ^{2}+\left\Vert \nabla g_{t}\wedge\nabla g_{t}\right\Vert ^{2}+\cdots+\left\Vert \nabla g_{t}\begin{array}{c}
n\\
\overbrace{\wedge\cdots\wedge}\end{array}\nabla g_{t}\right\Vert ^{2}}d\mathcal{L}^{n}\\
 & = & \int_{\Omega\backslash D_{\epsilon}}\left(\sum_{i}<\nabla_{i}g_{\epsilon},\nabla_{i}\eta(x/R)>+\sum_{k=2}^{n}\,\sum_{i_{1}<\cdots<i_{k},l,m}(-1)^{i_{l}+i_{m}}\left\langle \nabla_{i_{1}}g_{\epsilon}\wedge\begin{array}{c}
i_{l}\\
\cdots\\
\wedge\end{array}\wedge\nabla_{i_{k}}g_{\epsilon},\right.\right.\\
 &  & \left.\left.\left.,\nabla_{i_{1}}g_{\epsilon}\wedge\begin{array}{c}
i_{m}\\
\cdots\\
\wedge\end{array}\wedge\nabla_{i_{2}}g_{\epsilon}>\right\rangle <\nabla_{i_{l}}g_{\epsilon},\nabla_{i_{m}}\eta(x/R)>\right)\right/\\
 &  & \left(\sqrt{1+\left\Vert \nabla g_{\epsilon}\right\Vert ^{2}+\left\Vert \nabla g_{\epsilon}\wedge\nabla g_{\epsilon}\right\Vert ^{2}+\cdots+\left\Vert \nabla g_{\epsilon}\begin{array}{c}
n\\
\overbrace{\wedge\cdots\wedge}\end{array}\nabla g_{\epsilon}\right\Vert ^{2}}\right)d\mathcal{L}^{n}\end{eqnarray*}

Since the quadratic form $A$, at a fixed point $x$, defined by $(v,w)\mapsto<\nabla_{v}g_{\epsilon},\nabla_{w}g_{\epsilon}>:=<Av,w>$
is symmetric, there is an orthonormal basis $\{ e_{i}\}$ of $T_{*}(U,x)$
for which $A_{i}:=\nabla_{e_{i}}g_{\epsilon}=Ae_{i}$ are mutually
orthogonal, simplifying the calculations above somewhat. \begin{eqnarray*}
0 & = & \int_{U\backslash D_{\epsilon}}\left(\sum_{i}<A_{i},\nabla_{i}\eta(x/R)>+\right.\\
 &  & \left.\left.+\sum_{k=2}^{n}\,\sum_{i_{1}<\cdots<i_{k},l=1\ldots k}\left\Vert A_{i_{1}}\right\Vert ^{2}\begin{array}{c}
i_{l}\\
\cdots\\
\wedge\end{array}\left\Vert A_{i_{k}}\right\Vert ^{2}<A_{i_{l}},\nabla_{i_{l}}\eta((x-x_{0})/R)>\right)\right/\\
 &  & \left(\sqrt{1+\sum_{k=1}^{n}\,\sum_{i_{1}<\cdots<i_{k}}\left\Vert A_{i_{1}}\right\Vert ^{2}\cdots\left\Vert A_{i_{k}}\right\Vert ^{2}}\right)d\mathcal{L}^{n}\end{eqnarray*}

Additionally, since \[
1+\sum_{k=1}^{n}\,\sum_{i_{1}<\cdots<i_{k}}\left\Vert A_{i_{1}}\right\Vert ^{2}\cdots\left\Vert A_{i_{k}}\right\Vert ^{2}=\prod_{i=1}^{n}\left(1+\left\Vert A_{i}\right\Vert ^{2}\right),\]
and similarly, for each $j$ \[
1+\sum_{k=1}^{n-1}\,\sum_{i_{1}<\cdots<i_{k},i_{l}\neq j}\left\Vert A_{i_{1}}\right\Vert ^{2}\cdots\left\Vert A_{i_{k}}\right\Vert ^{2}=\prod_{i=1,i\neq j}^{n}\left(1+\left\Vert A_{i}\right\Vert ^{2}\right),\]
 as a functional over the base, \begin{eqnarray*}
0 & = & \int_{U\backslash D_{\epsilon}}\sum_{j=1}^{n}\frac{<A_{j},\nabla_{j}\eta(x/R)>}{1+\left\Vert A_{j}\right\Vert ^{2}}\left(\sqrt{1+\sum_{k=1}^{n}\,\sum_{i_{1}<\cdots<i_{k}}\left\Vert A_{i_{1}}\right\Vert ^{2}\cdots\left\Vert A_{i_{k}}\right\Vert ^{2}}\right)d\mathcal{L}^{n}\end{eqnarray*}
As a parametric integrand, the Euler-Lagrange equations simplify,
in this basis at each point, to \begin{eqnarray*}
0 & = & \left.\frac{d}{dt}\right|_{0}\int_{\pi^{-1}\left(U\backslash D_{\epsilon}\right)}f_{\epsilon}(\overrightarrow{T^{t}})d\left\Vert T^{t}\right\Vert \\
 & = & \int_{\pi^{-1}\left(U\backslash D_{\epsilon}\right)}\sum_{i}\frac{<\nabla_{i}g_{\epsilon},\nabla_{i}\eta(x/R)>}{1+\left\Vert \nabla_{i}g_{\epsilon}\right\Vert ^{2}}d\left\Vert T_{\epsilon}\right\Vert ,\end{eqnarray*}
where $g_{\epsilon}$ is the BV-carrier of $T_{\epsilon}$ on the
good set. Note that, although not explicitly manifest, the derivative
of $d\left\Vert T^{t}\right\Vert $ with respect to $t$ is included
in the formula above, since the preceding calculations are nonparametric
on the good set. 

On the bad set, the deformation is \begin{eqnarray*}
\left.\frac{d}{dt}\right|_{0}\int_{\pi^{-1}\left(D_{\epsilon}\right)}f_{\epsilon}(\overrightarrow{T^{t}})d\left\Vert T^{t}\right\Vert  & = & \int_{\pi^{-1}\left(D_{\epsilon}\right)}<\overrightarrow{T_{\epsilon}},\overrightarrow{h_{t,R}}>d\left\Vert T_{\epsilon}\right\Vert +B,\end{eqnarray*}
where, since the deformation is vertical, $B=0$. This follows since
a vertical deformation, that is, $T^{t}:=(H_{t})_{\# }(T)$ for $H_{t}(x,y)=(x,y+t\eta(x/R))$,
the boundary of the set where the penalty energy is nonzero, and the
penalty-energy $\mathcal{H}_{\epsilon}$ itself, will not change under
such a deformation. Also, the mass of that part of $T_{\epsilon}$
which is vertical (for which $\pi_{\#}(\overrightarrow{T_{\epsilon}})=0$)
will also remain unchanged under such a deformation.

\section{Squash-deformation }

Let $E$ be the cylindrical excess of the penalty-minimizer $T_{\epsilon}$,\[
E:=Exc(T_{\epsilon};R,x_{0}):=\frac{1}{R^{n}}\left(\mathcal{M}(T\rest\pi^{-1}(B(x_{0},R)))-\mathcal{M}(\pi_{\#}(T\rest\pi^{-1}(B(x_{0},R)))\right),\]
and for a given $R,\,0<R<1$, define the non-homothetic dilation $\phi_{R}(x,y)=(\frac{x}{R},\frac{y}{\sqrt{E}R})=(X,Y)$
of the cylinder $\pi^{-1}(B(x_{0},R))$ (we restrict to a coordinatizable
neighborhood, so that the fiber can be considered to be a compact
set within $\R^{j}$, and we assume without loss of generality that
$x_{0}=0$), and set $T_{\epsilon,R}:=\left(\phi_{R}\right)_{\#}\left(T_{\epsilon}\rest\pi^{-1}(B(x_{0},R))\right)$.
$T_{\epsilon,R}$ minimizes the penalty functional $\mathcal{F}_{\epsilon,R}$
defined by\begin{equation}
\mathcal{F}_{\epsilon,R}(S):=\int_{\pi^{-1}(B(x,R))}E^{-1}R^{-n}f_{\epsilon}\left(\overrightarrow{\left(\phi_{R}^{-1}\right)_{\#}(S)}\right)d\left\Vert \left(\phi_{R}^{-1}\right)_{\#}S\right\Vert ,\label{eq:squash-deformation}\end{equation}
which contracts the current $S$ back to the cylinder of radius $R$,
evaluates the original penalty functional there, and scales to compensate
for the factors of $R$ and some of the factors of $E$. Consider
the Euler-Lagrange equations of this functional, on $\widetilde{\Gamma}(B(x_{0},1)\times\R^{k})$.
Applying a vertical deformation as before,\begin{eqnarray*}
 &  & \frac{d}{dt}\mathcal{F}_{\epsilon R}\left(\left(h_{t}\right)_{\#}\left(T_{\epsilon,R}\right)\right)\\
 & = & \frac{d}{dt}\int_{\pi^{-1}(B(x_{0},R))}E^{-1}R^{-n}f_{\epsilon}\left(\overrightarrow{\left(\phi_{R}^{-1}\right)_{\#}\left(\left(h_{t}\right)_{\#}\left(T_{\epsilon,R}\right)\right)}\right)d\left\Vert \left(\phi_{R}^{-1}\right)_{\#}\left(h_{t}\right)_{\#}\left(T_{\epsilon,R}\right)\right\Vert \\
 & = & \frac{d}{dt}\int_{\pi^{-1}(B(x_{0},R))}E^{-1}R^{-n}f_{\epsilon}\left(\overrightarrow{\left(\left(h_{t,R,E}\right)_{\#}\left(T_{\epsilon}\right)\right)}\right)d\left\Vert \left(\left(h_{t,R}\right)_{\#}\left(T_{\epsilon}\right)\right)\right\Vert ,\end{eqnarray*}
where $h_{t,R,E}(x,y)=(x,y+\sqrt{E}Rt\eta(x/R))$

For a given $x$ in the good set, then for sufficiently small $R$
this integral consists of two pieces, the integral over the good set
within $B(x,R)$, which is the integral of a $C^{1}$ graph, and the
integral over the bad set, which shrinks with $\epsilon$. 

\begin{case}
On the good set, where $g_{t}$ and $G_{t}(X):=g_{t}(RX)/(\sqrt{E}R)$
are $C^{1}$, denote also by $\nabla_{i}G_{t}$ the covariant derivative
of $G_{t}$ in the direction of $\partial/\partial X_{i}$ on the
ball $B(0,1)$, with the metric stretched by the factor of $1/R$,
and similarly for other maps. For maps defined on $B(0,1)$, the notation
$\nabla_{i}$ will refer to covariant differentiation with respect
to $\partial/\partial X_{i}$, and for maps defined on $B(x_{0},R)$,
$\nabla_{i}$ will refer to covariant differentiation with respect
to $\partial/\partial x_{i}$. 
\end{case}
\begin{eqnarray*}
 &  & \frac{d}{dt}\left.\mathcal{F}_{\epsilon R}\right|_{B(x_{0},1)\backslash\phi_{R}(D_{\epsilon})}\left(\left(h_{t}\right)_{\#}\left(T_{\epsilon,R}\right)\right)\\
 & = & \frac{d}{dt}\int_{\pi^{-1}(B(x_{0},R)\backslash D_{\epsilon})}E^{-1}R^{-n}f_{\epsilon}\left(\overrightarrow{\left(\phi_{R}^{-1}\right)_{\#}\left(\left(h_{t}\right)_{\#}\left(T_{\epsilon,R}\right)\right)}\right)d\left\Vert \left(\phi_{R}^{-1}\right)_{\#}\left(h_{t}\right)_{\#}\left(T_{\epsilon,R}\right)\right\Vert \\
 & = & \int_{\pi^{-1}(B(x_{0},R)\backslash D_{\epsilon})}E^{-1}R^{-n}\frac{d}{dt}\left[f_{\epsilon}\left(\overrightarrow{\left(\left(h_{t,R,E}\right)_{\#}\left(T_{\epsilon}\right)\right)}\right)d\left\Vert \left(\left(h_{t,R,E}\right)_{\#}\left(T_{\epsilon}\right)\right)\right\Vert \right]\\
 & = & \int_{B(x_{0},R)\backslash D_{\epsilon}}E^{-1}R^{-n}\sum_{i}\frac{\sqrt{\Pi_{j}\left(1+\left\Vert \nabla_{j}g_{t}\right\Vert ^{2}\right)}}{1+\left\Vert \nabla_{i}g_{t}\right\Vert ^{2}}<\nabla_{i}g_{t},\nabla_{i}\eta_{R}>d\mathcal{L}^{n}\\
 & = & \int_{B(0,1)\backslash\phi_{R}(D_{\epsilon})}E^{-1}\sum_{i}\left.\frac{\sqrt{\Pi_{j}\left(1+\left\Vert \nabla_{j}g_{t}\right\Vert ^{2}\right)}}{1+\left\Vert \nabla_{i}g_{t}\right\Vert ^{2}}\right|_{x=RX}<\sqrt{E}\nabla_{i}G_{t},\sqrt{E}\nabla_{i}\eta>d\mathcal{L}^{n}\\
 & = & \int_{B(0,1)\backslash\phi_{R}(D_{\epsilon})}\sum_{i}\left.\frac{\sqrt{\Pi_{j}\left(1+\left\Vert \nabla_{j}g_{t}\right\Vert ^{2}\right)}}{1+\left\Vert \nabla_{i}g_{t}\right\Vert ^{2}}\right|_{x=RX}<\nabla_{i}G_{t},\nabla_{i}\eta>d\mathcal{L}^{n}.\end{eqnarray*}
Now, as $R\rightarrow0$, the integral formally becomes \begin{eqnarray*}
 & = & \int_{B(0,1)\backslash\lim_{R\rightarrow0}\phi_{R}(D_{\epsilon})}\sum_{i}\frac{\sqrt{\Pi_{j}\left(1+a_{j}^{2}\right)}}{1+A_{i}^{2}}<\nabla_{i}G_{t},\nabla_{i}\eta>d\mathcal{L}^{n},\end{eqnarray*}
where $a_{j}^{2}$ are the critical values of the quadratic form $(v,w)\mapsto<\nabla_{v}g_{\epsilon},\nabla_{w}g_{\epsilon}>:=<Av,w>$
as before, for unit vectors $v$ and $w$, defining a linear operator
$A$ as at the end of the previous section. $g_{\epsilon}$ is the
BV-carrier of the rectifiable section $T_{\epsilon}$. The operator
$\mathcal{A}=\sqrt{det(I+A)}(I+A)^{-1}$ will by elementary calculation
have the same eigenvectors as $A$, and eigenvalues: \begin{equation}
<\mathcal{A}e_{i},e_{i}>:=\frac{\sqrt{\Pi_{j}\left(1+a_{j}^{2}\right)}}{(1+a_{i}^{2})}.\label{eq:A_ij}\end{equation}

\begin{case}
On the bad set:
\end{case}
\begin{lem}
\label{lem:lemma 3a} Given $E=Exc(T_{\epsilon},R)$, and for any
$\epsilon>0$, \begin{eqnarray*}
\left\Vert T_{\epsilon}\rest\pi^{-1}(D_{\epsilon}\cap B(x_{0},R))\right\Vert  & \leq & ER^{n}+\epsilon A/\left|\log(\epsilon)\right|.\end{eqnarray*}

\end{lem}
\begin{pf}
\begin{eqnarray*}
\left\Vert T_{\epsilon}\rest\pi^{-1}(D_{\epsilon}\cap B(x_{0},R))\right\Vert  & = & \left(\left\Vert T_{\epsilon}\rest\pi^{-1}(D_{\epsilon}\cap B(x_{0},R))\right\Vert -\left\Vert D_{\epsilon}\cap B(x_{0},R)\right\Vert \right)\\
 &  & +\left\Vert D_{\epsilon}\cap B(x_{0},R)\right\Vert \\
 & \leq & ER^{n}+\left\Vert D_{\epsilon}\cap B(x_{0},R)\right\Vert \\
 & = & ER^{n}+\epsilon A/\left|\log(\epsilon)\right|,\end{eqnarray*}
where the first inequality follows from the fact that the excess is
that same difference between the mass of $T_{\epsilon}$ and its projection
(multiplicity 1) over a larger area than $D_{\epsilon}\cap B(x_{0},R)$. 
\hfill{}$\square$\end{pf}
Conversely, the excess $E$ will give a bound on the measure of $D_{\epsilon}$,
which will allow us to re-estimate the mass $\left\Vert T_{\epsilon}\rest\pi^{-1}(D_{\epsilon}\cap B(x_{0},R))\right\Vert $
in terms only of the excess. Recall that $A\ll B$ denotes that $A$
is bounded above by a constant times $B$, where that constant is
independent of the variables included in $A$ and $B$.

\begin{lem}
\label{lem:d-epsilon bound}$\left\Vert D_{\epsilon}\cap B(x_{0},R)\right\Vert \ll ER^{n}$.
\end{lem}
\begin{pf}
On the slightly smaller set $B_{\epsilon}\subset D_{\epsilon}$, $B_{\epsilon}:=\left\{ x\left|h_{\epsilon}((\overrightarrow{T_{\epsilon}})_{z})>0\,\textrm{for some\,}z\in\pi^{-1}(x)\right.\right\} $,
there will be at least 3 points in $\pi^{-1}(x)\cap Supp(T_{\epsilon})$
for \emph{a.e.} $x\in B_{\epsilon}$, because homologically $\pi_{\#}(T_{\epsilon})=1[B(x_{0},R)]$
and, where $h_{\epsilon}\neq0$, $\pi_{*}(\overrightarrow{T_{\epsilon}})=-1\overrightarrow{\R^{n}}$,
applying the constancy theorem. Thus, \begin{eqnarray*}
ER^{n} & \geq & \left\Vert T_{\epsilon}\rest\pi^{-1}(B_{\epsilon}\cap B(x_{0},R))\right\Vert -\left\Vert B_{\epsilon}\cap B(x_{0},R)\right\Vert \\
 & \geq & 2\left\Vert B_{\epsilon}\cap B(x_{0},R)\right\Vert ,\end{eqnarray*}
and the Lemma follows from the fact that $\left\Vert D_{\epsilon}\cap B(x_{0},R)\right\Vert \leq2\left\Vert B_{\epsilon}\cap B(x_{0},R)\right\Vert $,
by Proposition (\ref{pro:g}).
\hfill{}$\square$\end{pf}
\begin{rem}
On cursory examination, this Lemma would seem to imply that there
is a relationship between the excess and the penalty parameter $\epsilon$,
that is, the excess could not be chosen arbitrarily small unless $\epsilon$
is itself sufficiently small. Since, however, $D_{\epsilon}$ can
be empty independent of $\epsilon$, that is not necessarily the case.
\end{rem}
\begin{cor}
\label{lem:lemma 3a-improved} Given $\epsilon>0$ and $E=Exc(T_{\epsilon},R)$,
\begin{eqnarray*}
\left\Vert T_{\epsilon}\rest\pi^{-1}(D_{\epsilon}\cap B(x_{0},R))\right\Vert  & \ll & ER^{n}.\end{eqnarray*}

\end{cor}
\begin{pf}
If $\mathbf{e}$ is the unique unit horizontal $n$-plane so that
$\pi_{*}(\mathbf{e)}=*dV_{M}$\begin{eqnarray*}
ER^{n} & := & \left\Vert T_{\epsilon}\rest\pi^{-1}(B(x_{0},R))\right\Vert -\left\Vert B(x_{0},R)\right\Vert \\
 & = & \int_{\pi^{-1}(B(x_{0},R))}\left(1-<\overrightarrow{T_{\epsilon}},\mathbf{e}>\right)d\left\Vert T_{\epsilon}\right\Vert \\
 & \geq & \int_{\pi^{-1}(D_{\epsilon}\cap B(x_{0},R))}\left(1-<\overrightarrow{T_{\epsilon}},\mathbf{e}>\right)d\left\Vert T_{\epsilon}\right\Vert \\
 & = & \left\Vert T_{\epsilon}\rest\pi^{-1}(D_{\epsilon}\cap B(x_{0},R))\right\Vert -\left\Vert D_{\epsilon}\cap B(x_{0},R)\right\Vert \\
 & \gg & \left\Vert T_{\epsilon}\rest\pi^{-1}(D_{\epsilon}\cap B(x_{0},R))\right\Vert -ER^{n}\end{eqnarray*}
by Lemma (\ref{lem:d-epsilon bound}).
\hfill{}$\square$\end{pf}
In addition, we have

\begin{prop}
\begin{eqnarray*}
\frac{d}{dt}\left.\mathcal{F}_{\epsilon R}\right|_{\phi_{R}(D_{\epsilon})}\left(\left(h_{t}\right)_{\#}\left(T_{\epsilon,R}\right)\right) & \leq & C\sqrt{E}.\end{eqnarray*}

\end{prop}
\begin{pf}
\begin{eqnarray*}
 &  & \frac{d}{dt}\left.\mathcal{F}_{\epsilon R}\right|_{\phi_{R}(D_{\epsilon})}\left(\left(h_{t}\right)_{\#}\left(T_{\epsilon,R}\right)\right)\\
 & = & \frac{d}{dt}\int_{\pi^{-1}(D_{\epsilon})}E^{-1}R^{-n}f_{\epsilon}\left(\overrightarrow{\left(\phi_{R}^{-1}\right)_{\#}\left(\left(h_{t}\right)_{\#}\left(T_{\epsilon,R}\right)\right)}\right)d\left\Vert \left(\phi_{R}^{-1}\right)_{\#}\left(h_{t}\right)_{\#}\left(T_{\epsilon,R}\right)\right\Vert \\
 & = & \int_{\pi^{-1}(D_{\epsilon})}E^{-1}R^{-n}\frac{d}{dt}\left[f_{\epsilon}\left(\overrightarrow{\left(\left(h_{t,R,E}\right)_{\#}\left(T_{\epsilon}\right)\right)}\right)d\left\Vert \left(\left(h_{t,R,E}\right)_{\#}\left(T_{\epsilon}\right)\right)\right\Vert \right]\\
 & = & \int_{\pi^{-1}(D_{\epsilon})}E^{-1}R^{-n}<\overrightarrow{T_{\epsilon}},\overrightarrow{h_{t,R,E}}>d\left\Vert T_{\epsilon}\right\Vert \\
 & \leq & \int_{\pi^{-1}(D_{\epsilon})}E^{-1}R^{-n}\sqrt{E}d\left\Vert T_{\epsilon}\right\Vert \\
 & \leq & C\sqrt{E}.\end{eqnarray*}

\hfill{}$\square$\end{pf}

\section{Technical estimates}

There are a number of technical estimates we will need of higher Sobolev
and $L^{p}$ norms for the BV carrier $f$ of $T_{\epsilon}$ over
$B(x_{0},R)$. The notation is as in the previous section. These results
are all slight modifications of results in \cite{Bombieri}. The present
situation is, unfortunately, slightly different from that considered
by Bombieri, so that the statements, and proofs, need to be altered.

Following \cite{Bombieri}, first we show that 

\begin{lem}
\label{lem:lemma1}\[
\int_{B(x_{0},R)}\left(\left\Vert (dx,df)\right\Vert -1\right)d\mathcal{L}^{n}\leq ER^{n}.\]

\end{lem}
\begin{pf}
If $\eta$ is smooth and of compact support in the interior of $B(x_{0},R)$,
then \begin{eqnarray*}
\int\eta D_{i}f_{j} & = & -\int\frac{\partial\eta}{\partial x_{i}}f_{j}\\
 & = & -T_{\epsilon}\left(y_{j}\frac{\partial\eta}{\partial x_{i}}dx_{1}\wedge\cdots\wedge dx_{n}\right)\\
 & = & T_{\epsilon}\left((-1)^{i}y_{j}d\left(\eta dx_{1}\wedge\cdots\wedge\widehat{dx_{i}}\wedge\cdots\wedge dx_{n}\right)\right)\\
 & = & T_{\epsilon}\left((-1)^{i}d\left(y_{j}\eta dx_{1}\wedge\cdots\wedge\widehat{dx_{i}}\wedge\cdots\wedge dx_{n}\right)\right)+\\
 &  & +T_{\epsilon}\left((-1)^{i-1}\eta dy_{j}\wedge dx_{1}\wedge\cdots\wedge\widehat{dx_{i}}\wedge\cdots\wedge dx_{n}\right)\\
 & = & (-1)^{i}\partial T_{\epsilon}\left(y_{j}\eta dx_{1}\wedge\cdots\wedge\widehat{dx_{i}}\wedge\cdots\wedge dx_{n}\right)+\\
 &  & +T_{\epsilon}\left((-1)^{i-1}\eta dy_{j}\wedge dx_{1}\wedge\cdots\wedge\widehat{dx_{i}}\wedge\cdots\wedge dx_{n}\right)\\
 & = & T_{\epsilon}\left((-1)^{i-1}\eta dy_{j}\wedge dx_{1}\wedge\cdots\wedge\widehat{dx_{i}}\wedge\cdots\wedge dx_{n}\right).\end{eqnarray*}
Thus, by the definition of mass and the definition of $f$ as the
BV-carrier, \begin{equation}
\sup\int_{B(x_{0},R)}\left(\eta_{0}dx_{1}\wedge\cdots\wedge dx_{n}+\sum_{ij}\eta_{ij}D_{i}f_{j}\right)\leq M(T_{\epsilon}\rest\pi^{-1}(B(x_{0},R)))\label{eq:eta-0}\end{equation}
where the supremum is over all $(\eta_{0},\eta_{ij})$ of pointwise
norm less than or equal to 1. Since that supremum on the left is the
total variation of $(dx,df)$, subtracting $\int_{B(x_{0},R)}1d\mathcal{L}^{n}$
from both sides yields the statement.
\hfill{}$\square$\end{pf}
\begin{lem}
\label{lem:lemma2}\[
\int_{B(x_{0},R)}\left\Vert df\right\Vert \ll\sqrt{E}R^{n}.\]
Thus, there is a $y^{*}$ so that\[
\int_{B(x_{0},R)}\left|f(x)-y^{*}\right|d\mathcal{L}^{n}\ll\sqrt{E}R^{n+1}.\]

\end{lem}
\begin{pf}
In the inequality (\ref{eq:eta-0}), set $\eta_{0}=1-\tau$, $\tau>0$,
put all but the $D_{i}f_{j}$ terms on the right hand side, and we
get \[
\int_{B(x_{0},R)}\left(\sum_{ij}\eta_{ij}D_{i}f_{j}\right)\leq(\omega_{n}\tau+E)R^{n}\]
for all $\eta_{ij}$with $\sum\eta_{ij}^{2}\leq2\tau-\tau^{2}$, so
\[
\int_{B(x_{0},R)}\left\Vert df\right\Vert \leq\frac{(\omega_{n}\tau+E)R^{n}}{\sqrt{2\tau-\tau^{2}}}.\]
 Choose $\tau=E/(E+\omega_{n})$, then \[
\int_{B(x_{0},R)}\left\Vert df\right\Vert \leq\sqrt{E+\omega_{n}}\sqrt{E}R^{n}.\]
The second inequality follows from the first by a Poincar\'{e}-type
inequality for BV functions, proved by the standard contradiction
argument using the compactness theorem for BV functions. 
\hfill{}$\square$\end{pf}
\begin{rem}
Note that the implicit constant in the $\ll$ of the statement of
the Lemma is independent of $E$.
\end{rem}
\begin{lem}
\label{lem:lemma3c}For each $\epsilon>0$, the bad set $D_{\epsilon}$
can be chosen so that $\left\Vert \nabla g\right\Vert \ll1/\sqrt{E}$
on $B(x_{0},R)\backslash D_{\epsilon}$.
\end{lem}
\begin{pf}
By Lemma (\ref{lem:lemma2}), there is a constant $C$ so that $\int_{B(x_{0},R)}\left\Vert df\right\Vert d\mathcal{L}^{n}\leq C\sqrt{E}R^{n}$.
For each $A>0$, $\left\Vert \left\{ \left.x\in B(x_{0},R)\right|\left\Vert df(x)\right\Vert >A\right\} \right\Vert <C\sqrt{E}R^{n}/A$.
Given $1>\epsilon>0$, enlarge the bad set $D_{\epsilon}$ to also
include $\left\{ \left.x\in B(x_{0},R)\right|\left\Vert df(x)\right\Vert >1/\sqrt{E}\right\} $,
which will still keep the measure of the bad set $\left\Vert D_{\epsilon}\right\Vert \ll ER^{n}$.
\hfill{}$\square$\end{pf}
\begin{lem}
\label{lem:lemma5} For each penalty-minimizer $T_{\epsilon}$, there
is a $\gamma_{1}>$0 so that if the excess $E<\gamma_{1}$, we have
\[
Supp(T_{\epsilon}\rest\pi^{-1}(B(x_{0},R'))\subset\left\{ \left|y-y^{*}\right|\leq E^{\frac{1}{4n}}R\right\} ,\]
where $R'=(1-E^{1/4n})R$.
\end{lem}
\begin{pf}
Initially, we need some basic estimates. 

From Corollary {[}\ref{lem:lemma 3a-improved}{]}, $\left\Vert D_{\epsilon}\right\Vert <CE$.
For any given $v$, \begin{eqnarray*}
 &  & (vR)\cdot meas\left(B(x_{0},R)\cap\left\{ x\notin D_{\epsilon}\left|\left.\left|f(x)-y^{*}\right|\right.\right.>vR\right\} \right)\\
 & < & \int_{B(x_{0},R)}\left|f(x)-y^{*}\right|d\mathcal{L}^{n}\\
 & \ll & E^{1/2}R^{n+1},\end{eqnarray*}
by (\ref{lem:lemma2}) for the last inequality. Then, \[
meas\left(B(x_{0},R)\cap\left\{ x\notin D_{\epsilon}\left|\left|f(x)-y^{*}\right|\right.>vR\right\} \right)\ll\frac{1}{v}E^{1/2}R^{n}.\]
This implies that \begin{eqnarray*}
 &  & \left\Vert T_{\epsilon}\right\Vert \left\{ z=(x,y)\left|\left|y-y^{*}\right|\right.>vR\right\} \\
 & \ll & \left\Vert T_{\epsilon}\right\Vert \rest\pi^{-1}(D_{\epsilon})+meas\left(B(x_{0},R)\cap\left\{ x\notin D_{\epsilon}\left|\left.\left|f(x)-y^{*}\right|\right.\right.>vR\right\} \right)+ER^{n}\\
 & \ll & (2ER^{n}+\frac{1}{v}E^{1/2})R^{n}.\end{eqnarray*}
The proof of the Lemma now follows by a contradiction argument. Choose
$v=\frac{1}{2}E^{1/4n}$ and suppose there is a $z_{0}\in supp(T)$,
$z_{0}=(x_{0},y_{0})$, with $\left|y_{0}-y^{*}\right|>2vR$, and
with $\left|x_{0}\right|<(1-2v)R$ (without loss of generality we
can take $x_{0}=0$). Then, \[
\left\{ z\left|\left|z-z_{0}\right|\leq vR\right.\right\} \subset\left\{ z=(x,y)\left|\left.\left|y-y^{*}\right|>vR\right.\right.,\;\left|x-x_{0}\right|<R\right\} \]
and so the previous inequality implies \[
\mathcal{M}\left(T_{\epsilon}\rest\left\{ z=(x,y)\left|\left|z-z_{0}\right|\right.\leq vR\right\} \right)\ll(E+\frac{1}{v}E^{1/2})R^{n}\]
 Now, the monotonicity result Proposition(\ref{pro:monotonicity})
implies that for $\epsilon>0$ sufficiently small \[
(vR)^{n}\ll\mathcal{M}\left(T_{\epsilon}\rest\left\{ z\left|\left.\left|z-z_{0}\right|\right.\right.\leq vR\right\} \right),\]
stringing these inequalities together implies\[
\frac{1}{2^{n}}E^{1/4}R^{n}\ll(E+2E^{1/2-1/4n})R^{n}\ll(E+2E^{1/2-1/4n})R^{n}.\]
However, since the constant implied in the $\ll$ of this inequality
is again independent of $E$, for sufficiently small $E$ this inequality
will fail. Thus, there is a sufficiently small $E$, $E\leq\gamma_{1}$,
for which there is no such $z_{0}$; that is, for which the statement
of the Lemma will hold.
\hfill{}$\square$\end{pf}
\begin{lem}
\label{lemma 6}Set \[
\overline{y}\,:=\frac{1}{\left\Vert B(x_{0},R/2)\backslash\left(D_{\epsilon}\cap B(x_{0},R/2)\right)\right\Vert }\int_{B(x_{0},R/2)\backslash\left(D_{\epsilon}\cap B(x_{0},R/2)\right)}fd\mathcal{L}^{n}.\]
Then\[
\int_{\pi^{-1}(B(x_{0},R/2))}\left|y-\overline{y}\right|^{2}d\left\Vert T_{\epsilon}\right\Vert \ll E^{1+1/2n}R^{n+2}+\int_{B(x_{0},R/2)\backslash D_{\epsilon}}\left|f-\overline{y}\right|^{2}d\mathcal{L}^{n}.\]

\end{lem}
\begin{pf}
We have, from the proof of Lemma (\ref{lem:lemma5}) that, for any
$s>0$, \begin{eqnarray*}
 &  & \left\Vert T_{\epsilon}\right\Vert \left(\pi^{-1}(B(x_{0},R/2))\cap\left\{ \left|y-\overline{y}\right|>s\right\} \right)\\
 & \ll & ER^{n}\\
 &  & +\textrm{meas$\left(B(x_{0},R/2)\backslash\left(D_{\epsilon}\cap B(x_{0},R/2)\right)\cap\left\{ \left|y-\overline{y}\right|>s\right\} \right)$},\end{eqnarray*}
where the first term on the right-hand side is a bound on the mass
over the bad set $D_{\epsilon}\cap B(x_{0},R/2)$. Set $Y=\sup_{y\in Supp(T\rest\pi^{-1}(B(x_{0},R/2)))}\left|y-\overline{y}\right|$,
and we have\begin{eqnarray*}
 &  & \int_{\pi^{-1}(B(x_{0},R/2))}\left|y-\overline{y}\right|^{2}d\left\Vert T_{\epsilon}\right\Vert \\
 & = & 2\int_{0}^{Y}s\mathcal{M}\left(T\rest\pi^{-1}(B(x_{0},R/2))\cap\left\{ \left|y-\overline{y}\right|>s\right\} \right)ds\\
 & \ll & Y^{2}ER^{n}+\int_{B(x_{0},R/2)\backslash D_{\epsilon}}\left|f-\overline{y}\right|^{2}d\mathcal{L}^{n}.\end{eqnarray*}
Choose $\overline{x}$ with $\left|\overline{x}-x_{0}\right|\leq R/2$
and so that $(\overline{x},\overline{y})$ is in the convex closure
of $supp(T_{\epsilon}\rest\pi^{-1}(B(x_{0},R/2)))$ for $\epsilon$
sufficiently small so that the estimates in Lemma(\ref{lem:lemma5})
hold. That lemma then implies that \[
Y\leq\sup_{\pi^{-1}(B(x_{0},R/2)}\left|y-y^{*}\right|+\left|y^{*}-\overline{y}\right|\leq2E^{1/4n}R.\]
For $\epsilon>0$ sufficiently small, substituting this inequality
in above yields the Lemma.
\hfill{}$\square$\end{pf}
The following result, unlike the others of this section, is not merely
closely modeled upon the results of \cite{Bombieri}, it is precisely
as given in that paper. See \cite{Bombieri} for the proof, where
it is Lemma 7.

\begin{lem}
\label{lem:lemma7}Let $0<\theta\leq1,$ $1\leq p<\frac{n}{n-1}$.
there is a constant $\tau=\tau(\theta,p)$ such that if $A$ is a
measurable subset of $B(x_{0},R)$, if \[
meas(A)\geq\theta\, meas(B(x_{0},R)),\]
if $h\in BV(B(x_{0},R)$, and if either \[
\int_{A}hd\mathcal{L}^{n}=0\;\textrm{or$\;\int_{A}\textrm{sign}(h)\left|h\right|^{1/2}d\mathcal{L}^{n}=0,$ }\]
then\[
\left(R^{-n}\int_{B(x_{0},R)}\left|h\right|^{p}d\mathcal{L}^{n}\right)^{1/p}\leq\tau R^{1-n}\int_{B(x_{0},R)}\left|Dh\right|d\mathcal{L}^{n}.\]

\end{lem}

\begin{lem}
\label{lem:lemma 8}For\[
\overline{y}\,:=\frac{1}{\left\Vert B(x_{0},R/2)\backslash\left(D_{\epsilon}\cap B(x_{0},R/2)\right)\right\Vert }\int_{B(x_{0},R/2)\backslash\left(D_{\epsilon}\cap B(x_{0},R/2)\right)}fd\mathcal{L}^{n},\]
as in Lemma (\ref{lemma 6}), and if $E<1$, $1\leq p<\frac{n}{n-1}$,
we have

\[
\int_{B(x_{0},R/2)\backslash(B(x_{0},R/2)\cap D_{\epsilon})}\left|f-\overline{y}\right|^{2p}d\mathcal{L}^{n}\ll_{p}R^{n+2p}E^{p(1+1/2n)}.\]

\end{lem}
\begin{pf}
We may assume that $\overline{y}=0$. For $\phi=\phi(x)dx_{1}\wedge\cdots\wedge dx_{n}$
a horizontal form, define currents $V_{j}$ by \[
V_{j}(\phi):=T_{\epsilon}(y_{j}|y_{j}|\phi)\]
 and represent it by integration as \[
V_{j}(\phi)=\int_{B(x_{0},R)}h_{j}(x)\phi\]
with $h_{j}\in BV(B(x_{0},R)).$ By the definition of the good set
$B(x_{0},R/2)\backslash(B(x_{0},R/2)\cap D_{\epsilon})$, $h_{j}=f_{j}|f_{j}|$
on the good set. If $\psi=\sum_{i}\psi_{i}dx_{1}\wedge\cdots\wedge\widehat{dx_{i}}\wedge\cdots\wedge dx_{n}$
is smooth, with compact support in the interior of $B(x_{0},R),$
we have \begin{eqnarray*}
\partial V_{j}(\psi) & = & T_{\epsilon}(y_{j}|y_{j}|d\psi)\\
 & = & \partial T_{\epsilon}(y_{j}|y_{j}|\psi)-2T_{\epsilon}(|y_{j}|dy_{j}\wedge\psi)\\
 & = & -2T_{\epsilon}(|y_{j}|dy_{j}\wedge\psi).\end{eqnarray*}
If $\psi$ has compact support within $B(x_{0},R/2)$, \begin{eqnarray*}
\left|\partial V_{j}(\psi)\right| & \leq & 2\int_{B(x_{0},R/2)}\left|y_{j}\right|\sum_{i}\left|\psi_{i}\right|\left|\left\langle dy_{j}\wedge dx_{1}\wedge\cdots\wedge\widehat{dx_{i}}\wedge\cdots\wedge dx_{n},\overrightarrow{T_{\epsilon}}\right\rangle \right|d\left\Vert T_{\epsilon}\right\Vert \\
 & \leq & 2\left(\sup\left\Vert \psi\right\Vert \right)\int_{B(x_{0},R/2)}\left|y_{j}\right|\left(\sum_{i}\left|\left\langle dy_{j}\wedge dx_{1}\wedge\cdots\wedge\widehat{dx_{i}}\wedge\cdots\wedge dx_{n},\overrightarrow{T_{\epsilon}}\right\rangle \right|^{2}\right)^{1/2}d\left\Vert T_{\epsilon}\right\Vert \\
 & \leq & 2\left(\sup\left\Vert \psi\right\Vert \right)\left(\int_{C(x_{0},R/2)}\left|y_{j}\right|^{2}d\left\Vert T_{\epsilon}\right\Vert \right)^{1/2}\left(\int_{C(x_{0},R/2)}\left[1-\left|\left\langle dx,\overrightarrow{T_{\epsilon}}\right\rangle \right|^{2}\right]d\left\Vert T_{\epsilon}\right\Vert \right)^{1/2}\\
 & \leq & 2\left(\sup\left\Vert \psi\right\Vert \right)\left(\int_{C(x_{0},R/2)}\left|y_{j}\right|^{2}d\left\Vert T_{\epsilon}\right\Vert \right)^{1/2}\left(2ER^{n}\right)^{1/2}.\end{eqnarray*}
Lemma (\ref{lemma 6}) and this inequality implies that \begin{eqnarray*}
\int_{B(x,R/2)}\left|Dh_{j}\right|d\mathcal{L}^{n} & = & M\left(\left(\partial V_{j}\right)\rest B(x_{0},R/2)\right)\\
 & \ll & \left(\int_{C(x_{0},R/2)}\left|y_{j}\right|^{2}d\left\Vert T_{\epsilon}\right\Vert \right)^{1/2}\left(2ER^{n}\right)^{1/2}\\
 & \ll & \left(2ER^{n}\right)^{1/2}\left(E^{1+1/2n}R^{n+2}+\int_{B(x_{0},R/2)\backslash D_{\epsilon}}\left|h_{j}\right|d\mathcal{L}^{n}\right)^{1/2}.\end{eqnarray*}
Now apply Lemma(\ref{lem:lemma7}) with $A:=B(x_{0},R/2)\backslash D_{\epsilon}$
inside of $B(x_{0},R/2)$, which implies that \[
\int_{A}\left|h_{j}\right|d\mathcal{L}^{n}\ll R\int_{B(x_{0},R/2)}\left|Dh_{j}\right|d\mathcal{L}^{n}.\]
 Combining this with the previous inequality,\[
\left(\int_{B(x_{0},R/2)}\left|Dh_{j}\right|d\mathcal{L}^{n}\right)^{2}\ll\left(2ER^{n}\right)\left(E^{1+1/2n}+R\int_{B(x_{0},R/2)}\left|Dh_{j}\right|d\mathcal{L}^{n}\right),\]
which by the quadratic formula and the fact that $E<1$ implies that
\[
\int_{B(x_{0},R/2)}\left|Dh_{j}\right|d\mathcal{L}^{n}\ll E^{1+1/2n}R^{n+1}.\]
Applying Lemma (\ref{lem:lemma7}) gives the statement.
\hfill{}$\square$\end{pf}
\begin{lem}
\label{lemma_10} There is an $r$ with $R/4\leq r\leq R/3$, for
which, given $0<\mu\leq1$, there is a current $S$ so that
\end{lem}
\begin{enumerate}
\item $\partial(S\rest C(x_{0},r))=\partial(T_{\epsilon}\rest C(x_{0},r))$,
\item $\partial(\pi_{\#}(S\rest C(x_{0},R))=\partial B(x_{0},R)$,
\item $diam(Supp(S\rest C(x_{0},R))\cup Supp(T_{\epsilon}\rest C(x_{0},r)))\leq R$,
\item $Exc(S,R)\ll\mu E+E^{1+1/2n}/\mu+\left.\int_{B(x_{0},R/2)\backslash D_{\epsilon}}\left|f-\overline{y}\right|^{2}d\mathcal{L}^{n}\right/(\mu R^{n+2})$.
\end{enumerate}
\begin{pf}
As before, normalize so that $\overline{y}=0$. If $S$ is any normal
current in $\Omega\times\R^{k}$, the slice \[
<S,r>:=\partial(S\rest C(x_{0},r))-(\partial S)\rest C(x_{0},r)\]
satisfies, for smooth functions $g$, \[
<S,r>\rest g=<S\rest g,r>\]
for almost every $r$, where $S\rest g(\phi):=S(g\phi)$, and \[
\int_{0}^{p}\mathcal{M}(<S,r>)dr\leq\mathcal{M}(S\rest C(x_{0},p))\]
 (cf. Morgan, p. 55). Applied to $T_{\epsilon}$, with $g=|y|^{2}$,
and $p=R/2$, we have\begin{eqnarray}
\int_{0}^{R/2}\left(\mathcal{M}(<T_{\epsilon},r>)-n\alpha_{n}r^{n-1}\right)dr & \leq & \mathcal{M}(T_{\epsilon}\rest C(x_{0},R/2))-\alpha_{n}R^{n}/2^{n}\nonumber \\
 & \le & \left(\frac{R}{2}\right)^{n}Exc(T_{\epsilon},R/2)\nonumber \\
 & \leq & R^{n}E,\label{eq:slice-1}\end{eqnarray}
and, from Lemma (\ref{lemma 6})\begin{eqnarray}
\int_{0}^{R/2}\mathcal{M}(<T_{\epsilon},r>\rest|y|^{2})dr & = & \int_{0}^{R/2}\mathcal{M}(<T_{\epsilon}\rest|y|^{2},r>)dr\nonumber \\
 & \leq & \int_{C(x_{0},R/2)}|y|^{2}d\left\Vert T_{\epsilon}\right\Vert \nonumber \\
 & \leq & E^{1+1/2n}R^{n+2}+\int_{B(x_{0},R/2)\backslash D_{\epsilon}}\left|f\right|^{2}d\mathcal{L}^{n}.\label{eq:slice 2}\end{eqnarray}
Note also that \[
\mathcal{M}(<T_{\epsilon},r>)-n\alpha_{n}r^{n-1}\geq\mathcal{M}(<\pi_{\#}(T_{\epsilon}),r>)-n\alpha_{n}r^{n-1}=0,\]
since $\pi_{\#}$ is mass-decreasing. 

From (\ref{eq:slice-1}), there is some $r$ with \[
\mathcal{M}(<T_{\epsilon},r>)-n\alpha_{n}r^{n-1}\ll ER^{n-1},\]
and due to the implicit constant in the inequality, such an $r$ can
be found in $[R/4,R/3]$. We can also find, using (\ref{eq:slice 2}),
a choice of $r\in[R/4,R/3]$ also satisfying \[
\mathcal{M}(<T_{\epsilon},r>\rest|y|^{2})\ll E^{1+1/2n}R^{n+1}+\frac{1}{R}\int_{B(x_{0},R/2)\backslash D_{\epsilon}}\left|f\right|^{2}d\mathcal{L}^{n}.\]

We now construct a comparison current. Set $S$ to be the current\[
S:=B(x_{0},(1-\mu)r)\times\{0\}+h_{\#}([1-\mu,1+\mu]\times<T_{\epsilon},r>)+(B(x_{0},R)-B(x_{0},(1+\mu)r))\times\{0\},\]
where \[
h(t,x,y)=(tx,y-|t-1|y/\mu).\]
$S$ is a deformation of the horizontal current $B(x_{0},R)\times\{0\}$
that matches with the slice of $T_{\epsilon}$ at radius $r$, but
which is still flat off of an annulus of width $2\mu$. It is clear
from the construction that this current satisfies (1) and (2) of the
statement.

Since $\left|\partial h/\partial t\right|\leq(r^{2}+|y|^{2}/\mu^{2})^{1/2}$
(also cf. \cite[4.1.9]{GMT})\[
\mathcal{M}(h_{\#}([1-\mu,1+\mu]\times<T_{\epsilon},r>))\leq\int_{1-\mu}^{1+\mu}t^{n-1}\int(r^{2}+|y|^{2}/\mu^{2})^{1/2}d\left\Vert <T_{\epsilon},r>\right\Vert dt.\]
Performing the indicated integration with respect to $t$ and noting
that $(r^{2}+|y|^{2}/\mu^{2})^{1/2}\leq\left(r+\frac{|y|^{2}}{\mu^{2}r}\right)$,
\begin{eqnarray}
\mathcal{M}(h_{\#}([1-\mu,1+\mu]\times<T_{\epsilon},r>)) & \leq & \left(\frac{(1+\mu)^{n}-(1-\mu)^{n}}{n}\right)\int\left(r+\frac{|y|^{2}}{\mu^{2}r}\right)d\left\Vert <T_{\epsilon},r>\right\Vert \nonumber \\
 & \leq & \left(2n\mu\right)\left(r\mathcal{M}(<T_{\epsilon},r>)+\frac{1}{\mu^{2}r}\mathcal{M}(<T_{\epsilon},r>\rest|y|^{2})\right).\label{eq:comparison}\end{eqnarray}
Now, \begin{eqnarray*}
Exc(S,R) & = & \mathcal{M}(S\rest C(x_{0},R))/R^{n}-\alpha_{n}\\
 & = & \mathcal{M}(h_{\#}([1-\mu,1+\mu]\times<T_{\epsilon},r>))/R^{n}+\alpha_{n}((1-\mu)^{n}r^{n}+(1+\mu)^{n}r^{n})/R^{n}\\
 & \leq & 2n\mu\left(r\mathcal{M}(<T_{\epsilon},r>)+\frac{1}{\mu^{2}r}\mathcal{M}(<T_{\epsilon},r>\rest|y|^{2})\right)/R^{n}+\alpha_{n}(-2n\mu r^{n})/R^{n}\\
 & \ll & \mu\left(r(ER^{n-1})+\frac{1}{\mu^{2}r}\left(E^{1+1/2n}R^{n+1}+\frac{1}{R}\int_{B(x_{0},R/2)\backslash D_{\epsilon}}\left|f-\overline{y}\right|^{2}d\mathcal{L}^{n}\right)\right)/R^{n}\\
 & \ll & \mu E+\frac{1}{\mu}E^{1+1/2n}+\frac{1}{\mu R^{n+2}}\left(\int_{B(x_{0},R/2)\backslash D_{\epsilon}}\left|f-\overline{y}\right|^{2}d\mathcal{L}^{n}\right),\end{eqnarray*}
which is part (4) of the Lemma. 

Part (3) of the Lemma follows from Lemma (\ref{lem:lemma5}).
\hfill{}$\square$\end{pf}
\begin{lem}
\label{lemma_11}If $R\leq\gamma_{3}$, then, for $0<\mu\leq1$ chosen
as before, and if $E\leq\min\{\gamma_{1},(2/3)^{4n}\}$,\begin{eqnarray*}
Exc(T_{\epsilon},R/4) & \ll & \mu E(1+\frac{1}{2\epsilon})+E\left(\frac{E^{1/2n}}{\mu}(1+\frac{1}{2\epsilon})\right)+\left.\int_{B(x_{0},R/2)\backslash D_{\epsilon}}\left|f-\overline{y}\right|^{2}d\mathcal{L}^{n}\right/(\mu R^{n+2}).\end{eqnarray*}

\end{lem}
\begin{pf}
Again, suppose that $\overline{y}=0$. Let $S$ be as in Lemma (\ref{lemma_10}).
and set \[
\widetilde{T}:=T_{\epsilon}\rest C(x_{0},r)+S-S\rest C(x_{0},r),\]
which replaces $T_{\epsilon}$ by $S$ outside of the cylinder of
radius $r$, without introducing any interior boundaries by the construction
of $S$. Note that $\partial\widetilde{T}=\partial B(x_{0},R)\times\{0\}$.
By construction, monotonicity of the unnormalized excess, and the
choice of $r$, $R/4\leq r\leq R/3$, \[
(R/4)^{n}Exc(T_{\epsilon},R/4)\leq r^{n}Exc(T_{\epsilon},r)=r^{n}Exc(\widetilde{T},r)\leq R^{n}Exc(\widetilde{T},R).\]
By the definition of the penalty functional, \begin{eqnarray*}
Exc(\widetilde{T},R) & := & \left(\mathcal{M}(\widetilde{T})-\mathcal{M}(B(x_{0},R)\times\{0\})\right)/R^{n}\\
 & \leq & \left(\mathcal{F}_{\epsilon}(\widetilde{T})-\mathcal{F}_{\epsilon}(B(x_{0},R)\times\{0\})\right)/R^{n}.\end{eqnarray*}
Using minimality, \[
\mathcal{F}_{\epsilon}(T_{\epsilon}\rest C(x_{0},r))\leq\mathcal{F}_{\epsilon}(S\rest C(x_{0},r)),\]
so that \begin{eqnarray*}
\mathcal{F}_{\epsilon}(\widetilde{T}) & = & \mathcal{F}_{\epsilon}(T_{\epsilon}\rest C(x_{0},r)+(S-S\rest C(x_{0},r)))\\
 & \leq & \mathcal{F}_{\epsilon}(S).\end{eqnarray*}
Thus, \begin{eqnarray*}
Exc(\widetilde{T},R) & \leq & \left(\mathcal{F}_{\epsilon}(\widetilde{T})-\mathcal{F}_{\epsilon}(B(x_{0},R)\times\{0\})\right)/R^{n}\\
 & \leq & \left(\mathcal{F}_{\epsilon}(S)-\mathcal{M}(B(x_{0},R)\times\{0\})\right)/R^{n}.\end{eqnarray*}

Now, \begin{eqnarray*}
\mathcal{F}_{\epsilon}(S) & = & \mathcal{M}(B(x_{0},(1-\mu)r)\times\{0\})+\mathcal{M}((B(x_{0},R)-B(x_{0},(1+\mu)r))\times\{0\})\\
 &  & +\mathcal{F_{\epsilon}}(h_{\#}([1-\mu,1+\mu]\times<T_{\epsilon},r>)),\end{eqnarray*}
and the slice $<T_{\epsilon},r>$ is the graph of the $C^{1}$ function
$f$ on $\partial B(x_{0},r)\backslash(D_{\epsilon}\cap\partial B(x_{0},r))$.
The integral over the bad set $D_{\epsilon}\cap\partial B(x_{0},r)$
will, for some $r\in[R/4,R/3]$ consistent with all previous choices
of $r$, be bounded by the mass over that set plus $(12/R)(ER^{n})(\frac{1}{2\epsilon})=12R^{n-1}(E)(\frac{1}{2\epsilon})$
by Corollary (\ref{lem:lemma 3a-improved}) and the definition of
$\mathcal{F}_{\epsilon}$. So, similarly to equation (\ref{eq:comparison})the
proof of Lemma (\ref{lemma_10}), but using the height bound of Lemma
(\ref{lem:lemma5}) to bound $|y|$, along with the estimate for $|y|$
from Lemma (\ref{lemma 6}), \[
Supp(T_{\epsilon}\rest\pi^{-1}(B(x_{0},R'))\subset\left\{ \left|y-y^{*}\right|\leq E^{\frac{1}{4n}}R\right\} ,\]
 and \[
\sup_{\pi^{-1}(B(x_{0},R/2)}\left|y-y^{*}\right|+\left|y^{*}-\overline{y}\right|\leq2E^{1/4n}R.,\]
with $\overline{y}=0$, implying that $|y|<2E^{1/4n}R$, to bound
the contribution from the sloped sides of $S$ on the bad set,\begin{eqnarray*}
\mathcal{F_{\epsilon}}(h_{\#}([1-\mu,1+\mu]\times<T_{\epsilon},r>)) & \leq & \int_{B(x_{0},r(1+\mu))\backslash B(x_{0},r(1-\mu))}\mathcal{M}(h_{\#}([1-\mu,1+\mu]\times<T_{\epsilon},r>))\\
 & + & (r^{2}+|y|^{2}/\mu^{2})^{1/2}\left(\frac{(1+\mu)^{n}-(1-\mu)^{n}}{n}\right)(12R^{n-1})(E)(\frac{1}{2\epsilon})\\
 & \leq & \int_{B(x_{0},r(1+\mu))\backslash B(x_{0},r(1-\mu))}\mathcal{M}(h_{\#}([1-\mu,1+\mu]\times<T_{\epsilon},r>))\\
 &  & +2n\mu\left(r+\frac{2R^{2}E^{1/2n}}{\mu^{2}r}\right)(12R^{n-1})(E)(\frac{1}{2\epsilon}).\end{eqnarray*}
Combining this inequality with Lemma (\ref{lemma_10}), \begin{eqnarray*}
Exc(T_{\epsilon},R/4) & \leq & 4^{n}Exc(\widetilde{T},R)\\
 & \leq & 4^{n}\left(\mathcal{F}_{\epsilon}(S)-\mathcal{M}(B(x_{0},R)\times\{0\})\right)/R^{n}\\
 & \leq & 4^{n}\left(\mathcal{M}(S)-\mathcal{M}(B(x_{0},R)\times\{0\})\right.\\
 &  & \left.+2n\mu\left(r+\frac{2R^{2}E^{1/2n}}{\mu^{2}r}\right)(12R^{n-1})(E)(\frac{1}{2\epsilon})\right)/R^{n}\\
 & \leq & 4^{n}\left(Exc(S)+2n\mu\left(r+\frac{2R^{2}E^{1/2n}}{\mu^{2}r}\right)(12R^{-1})(E)(\frac{1}{2\epsilon})\right)\\
 & \ll & \mu E(1+\frac{1}{2\epsilon})+E\left(\frac{E^{1/2n}}{\mu}(1+\frac{1}{2\epsilon})\right)+\left.\int_{B(x_{0},R/2)\backslash D_{\epsilon}}\left|f-\overline{y}\right|^{2}d\mathcal{L}^{n}\right/(\mu R^{n+2}),\end{eqnarray*}
as required.
\hfill{}$\square$\end{pf}

\section{First variation of $\mathcal{F}_{\epsilon}(T)$}

Consider the deformations $(h_{t})_{\#}(T_{\epsilon})$ of $T_{\epsilon}$,
where $h_{t}$ is given by\[
h_{t}(x,y):=(x,y+t\sqrt{E}R\eta(x/R)),\]
for $-1<t<1,$ and $\eta$ smooth with compact support in $|X|<1$,
with $\left\Vert \nabla\eta\right\Vert \leq\beta$. Given the blow-up
map \[
\phi_{R}(x,y)=\left(\frac{x}{R},\frac{y}{\sqrt{E}R}\right):=(X,Y),\]
define $F:B(x_{0},1)\rightarrow\R^{j}$ by \[
F(X)=f(RX)/(\sqrt{E}R),\]
where $f$ is, as before, the BV-carrier of $T_{\epsilon}$. On the
good set, moreover, $G_{\epsilon}(X)=g_{\epsilon}(RX)/\sqrt{E}R,$
and so $\nabla_{X}G_{\epsilon}=\frac{1}{\sqrt{E}}\nabla_{x}g_{\epsilon}$,
where $g_{\epsilon}$ is the graph representing $T_{\epsilon}$ on
the good set.

\begin{lem}
\label{lemma 12} If $\eta(X)$ is smooth with compact support in
$|X|<1$, $|\nabla\eta|\leq1$, then given a deformation $h_{t}$
given by \[
h_{t}(x,y):=(x,y+t\sqrt{E}R\eta(x/R))\]
 and if $T_{\epsilon,R}=(\phi_{R})_{\#}(T_{\epsilon})$, where $\phi_{R}(x,y)=(x/R,y/(\sqrt{E}R))$,
then \[
\left|\frac{d}{dt}\mathcal{F}_{\epsilon,R}\left(\left(h_{t}\right)_{\#}\left(T_{\epsilon,R}\right)\right)-\int_{B(X_{0},1)}\sum_{i,k}\mathcal{A}_{ik}\left\langle \nabla_{i}\eta,\nabla_{k}F+t\nabla_{k}\eta\right\rangle d\mathcal{L}^{n}\right|\ll\sqrt{E}.\]

\end{lem}
\begin{pf}
By Lemma (\ref{lem:lemma 3a-improved}), and the definition of the
bad set $D_{\epsilon}$ in Proposition (\ref{pro:g}), we find a $C^{1}$
function $g_{\epsilon}:B(x_{0},R)\backslash D_{\epsilon}\rightarrow F$
whose graph agrees with $T_{\epsilon}$ over $B(x_{0},R)\backslash D_{\epsilon}$,
and $g_{t}(x):=g_{\epsilon}(x)+t\sqrt{E}R\eta(x/R)$. Then \begin{eqnarray*}
L(t) & := & \mathcal{F}_{\epsilon}((h_{t})_{\#}(graph(g_{\epsilon}))\rest(C(x_{0},R)\backslash\pi^{-1}(D_{\epsilon}))/(ER^{n}),\\
K(t) & := & \mathcal{F}_{\epsilon}((h_{t})_{\#}(T_{\epsilon})\rest(C(x_{0},R)\cap\pi^{-1}(D_{\epsilon}))/(ER^{n})\end{eqnarray*}
so that \[
\mathcal{F}_{\epsilon}((h_{t})_{\#}(T_{\epsilon}))/(ER^{n})=L(t)+K(t).\]

Apply the squash-deformation $\phi_{R}(x,y):=(x/R,y/(\sqrt{E}R))$.
If $T_{\epsilon,R}:=(\phi_{R})_{\#}(T_{\epsilon}\rest C(x_{0},R))$,
it will minimize the functional $\mathcal{F}_{\epsilon,R}$ defined
by \[
\mathcal{F}_{\epsilon,R}(S):=\mathcal{F}_{\epsilon}((\phi_{R}^{-1})_{\#}(S))/(ER^{n}),\]
 so that, on $T_{\epsilon,R}$, $\mathcal{F}_{\epsilon,R}(T_{\epsilon,R}):=\mathcal{F}_{\epsilon}(T_{\epsilon}\rest C(x_{0},R))/(ER^{n})$.
Explicitly, for $S$ a graph on $\pi^{-1}(\Omega)\subset C(X_{0},1)$,
\begin{eqnarray*}
\mathcal{F}_{\epsilon,R}(S) & = & \left\Vert (\phi_{R}^{-1})_{\#}(S)\right\Vert /(ER^{n})+\frac{1}{\epsilon ER^{n}}\mathcal{H}_{0}(S),\end{eqnarray*}
where $\mathcal{H}_{0}$ is as defined in the beginning of \S4.

On the good set, since the penalty term vanishes there,\begin{eqnarray*}
\frac{d}{dt}L(t) & = & \frac{d}{dt}\mathcal{F}_{\epsilon}((h_{t})_{\#}(graph(g_{\epsilon}))\rest(C(x_{0},R)\backslash\pi^{-1}(D_{\epsilon}))/(ER^{n})\\
 & = & \frac{d}{dt}\mathcal{M}((h_{t})_{\#}(graph(g_{\epsilon}))\rest(C(x_{0},R)\backslash\pi^{-1}(D_{\epsilon}))/(ER^{n})\\
 & = & \left.\int_{\pi^{-1}\left(B(x_{0},R)\backslash D_{\epsilon}\right)}\sum_{i}\frac{<\nabla_{i}g_{t},\nabla_{i}h>}{1+\left\Vert \nabla_{i}g_{t}\right\Vert ^{2}}d\left\Vert T_{t}\right\Vert \right/(ER^{n})\end{eqnarray*}
by {[}\ref{sub:Euler-Lagrange-equations-for-T-epsilon}{]}. Since
$g_{t}(x):=g_{\epsilon}(x)+t\sqrt{E}R\eta(x/R)$ and $h(x)=\frac{dg_{t}}{dt}=\sqrt{E}R\eta(x/R)$,
\begin{eqnarray*}
\frac{d}{dt}L(t) & = & \int_{\pi^{-1}\left(B(x_{0},R)\backslash D_{\epsilon}\right)}\sqrt{E}\sum_{i}\frac{<\nabla_{i}g_{\epsilon},\nabla_{i}\eta>+t\sqrt{E}<\nabla_{i}\eta,\nabla_{i}\eta>}{1+\left\Vert \nabla_{i}g_{t}\right\Vert ^{2}}d\left\Vert T_{t}\right\Vert /ER^{n}\\
 & = & \int_{B(x_{0},R)\backslash D_{\epsilon}}\frac{1}{\sqrt{E}R^{n}}\sum_{i}\frac{<\nabla_{i}g_{\epsilon},\nabla_{i}\eta>+t\sqrt{E}<\nabla_{i}\eta,\nabla_{i}\eta>}{1+\left\Vert \nabla_{i}g_{\epsilon}\right\Vert ^{2}+2t\sqrt{E}<\nabla_{i}g_{\epsilon},\nabla_{i}\eta>+t^{2}E\left\Vert \nabla_{i}\eta\right\Vert ^{2}}\\
 &  & \left\Vert (e_{1}+\nabla_{1}g_{\epsilon}+t\sqrt{E}\nabla_{1}\eta)\wedge\cdots\wedge(e_{n}+\nabla_{n}g_{\epsilon}+t\sqrt{E}\nabla_{n}\eta)\right\Vert d\mathcal{L}^{n}.\end{eqnarray*}
Now apply the squash-deformation $\phi_{R}(x,y)=(X,Y):=(x/R,y/(\sqrt{E}R))$.
Explicitly, for $S$ a graph, $S=graph(P(X))$ on $C(X_{0},1)$, \begin{eqnarray*}
\mathcal{F}_{\epsilon,R}(S) & = & \left\Vert (\phi_{R}^{-1})_{\#}(S)\right\Vert /(ER^{n})+\frac{1}{\epsilon ER^{n}}H_{0}(S)\\
 & = & \frac{1}{L}\int_{B(x_{0},R)}\sqrt{1+E\left\Vert \nabla_{i}P\right\Vert ^{2}+\cdots+E^{n}\left\Vert \nabla_{i_{1}}P\wedge\cdots\wedge\nabla_{i_{n}}P\right\Vert ^{2}}d\mathcal{L}^{n},\end{eqnarray*}
keeping in mind that the penalty term vanishes on graphs. Use a coordinate
system $\{ x^{1},\ldots,x^{n}\}$ so that the quadratic form $A_{\epsilon}(v,w)\mapsto\left.<\nabla_{v}g_{\epsilon},\nabla_{w}g_{\epsilon}>\right|_{x_{0}}$
is diagonalized, with eigenvalues $a_{i}^{2}$. The operator $\mathcal{A}_{\epsilon}:=\sqrt{det(I+A_{\epsilon})}(I+A_{\epsilon})^{-1}$
with the same eigenvectors but with eigenvalues $\mathcal{A}_{\epsilon,i}=\frac{\sqrt{\Pi_{j}(1+a_{j}^{2})}}{1+a_{i}^{2}}$
is the first term in the expansion of the previous expression.\begin{eqnarray*}
\frac{d}{dt}L(t) & = & \frac{d}{dt}\mathcal{F}_{\epsilon}((h_{t})_{\#}(graph(g_{\epsilon}))\rest(C(x_{0},R)\backslash\pi^{-1}(D_{\epsilon}))/(ER^{n})\\
 & =\\
 & = & \frac{d}{dt}\mathcal{F}_{\epsilon,R}((\phi_{R})_{\#}(h_{t})_{\#}(T_{\epsilon})\rest C(x_{0},R)\backslash\pi^{-1}(D_{\epsilon}))\\
 & = & \frac{d}{dt}\mathcal{F}_{\epsilon,R}\left(graph(G_{\epsilon}+t\eta)\rest\left(C(X_{0},1)\backslash\phi_{R}(\pi^{-1}(D_{\epsilon}))\right)\right)\\
 & = & \frac{1}{E}\int_{B(X_{0},1)\backslash\phi_{R}(D_{\epsilon})}\sum_{i}\frac{E<\nabla_{i}G_{\epsilon},\nabla_{i}\eta>+tE<\nabla_{i}\eta,\nabla_{i}\eta>}{1+E\left\Vert \nabla_{i}(G_{\epsilon}+t\eta)\right\Vert ^{2}}\cdot\\
 &  & \cdot\sqrt{1+E\left\Vert \nabla(G_{\epsilon}+t\eta)\right\Vert ^{2}+\cdots+E^{n}\left\Vert \nabla(G_{\epsilon}+t\eta)\wedge\cdots\wedge\nabla(G_{\epsilon}+t\eta)\right\Vert ^{2}}d\mathcal{L}^{n}\\
 & = & \int_{B(X_{0},1)\backslash\phi_{R}(D_{\epsilon})}\sum_{i}\left(<\nabla_{i}G_{\epsilon},\nabla_{i}\eta>+t<\nabla_{i}\eta,\nabla_{i}\eta>\right)\mathcal{A}_{ii}d\mathcal{L}^{n}+Q,\end{eqnarray*}
where the coordinate basis $\{ X_{1},\ldots,X_{n}\}$ is chosen at
each point to be an orthonormal eigenbasis of $(V,W)\mapsto<\nabla_{V}(G_{\epsilon}+t\eta),\nabla_{W}(G_{\epsilon}+t\eta)>$
and, at each point, $\nabla_{i}:=\nabla_{\partial/\partial X_{i} }$.
Since $\left\{ \nabla_{j}(G_{\epsilon}+t\eta)\right\} $ is orthogonal
by choice of basis, \begin{eqnarray*}
 &  & \sqrt{1+E\left\Vert \nabla(G_{\epsilon}+t\eta)\right\Vert ^{2}+\cdots+E\left\Vert \nabla(G_{\epsilon}+t\eta)\wedge\cdots\wedge\nabla(G_{\epsilon}+t\eta)\right\Vert ^{2}}\\
 & = & \sqrt{\Pi_{j}\left(1+E\left\Vert \nabla_{i}(G_{\epsilon}+t\eta)\right\Vert ^{2}\right)}\end{eqnarray*}
Choose $\{ V_{i}\}$ to be an eigenbasis of $\mathcal{A}$ as above,
that is, an eigenbasis of $(V,W)\mapsto<\nabla_{V}(G_{\epsilon}+0\eta),\nabla_{W}(G_{\epsilon}+0\eta)>$
at $X_{0}$. 

$Q$ is given simply as \begin{eqnarray*}
Q & := & \int_{B(X_{0},1)\backslash\phi_{R}(D_{\epsilon})}\sum_{i}\frac{<\nabla_{i}G_{\epsilon},\nabla_{i}\eta>+t<\nabla_{i}\eta,\nabla_{i}\eta>}{1+E\left\Vert \nabla_{i}(G_{\epsilon}+t\eta)\right\Vert ^{2}}\sqrt{\Pi_{j}\left(1+E\left\Vert \nabla_{j}(G_{\epsilon}+t\eta)\right\Vert ^{2}\right)}\\
 &  & -\left(<\nabla_{i}G_{\epsilon},\nabla_{i}\eta>+t<\nabla_{i}\eta,\nabla_{i}\eta>\right)\mathcal{A}_{ii}d\mathcal{L}^{n}\\
 & := & \int_{B(X_{0},1)\backslash\phi_{R}(D_{\epsilon})}\sum_{i}\left(<\nabla_{i}G_{\epsilon},\nabla_{i}\eta>+t<\nabla_{i}\eta,\nabla_{i}\eta>\right)Q_{i}d\mathcal{L}^{n}.\end{eqnarray*}
If \begin{eqnarray*}
Q_{i}(P_{1},\ldots,P_{n}) & := & \frac{\sqrt{\Pi_{j\neq i}\left(1+E\left\Vert P_{j}\right\Vert ^{2}\right)}}{\sqrt{1+E\left\Vert P_{i}\right\Vert ^{2}}}-\mathcal{A}_{ii},\end{eqnarray*}
$Q_{i}:=Q_{i}(\nabla_{1}(G_{\epsilon}+t\eta),\ldots,\nabla_{n}(G_{\epsilon}+t\eta))$,
then by a simple application of the mean value theorem at each $x$,
there is a $c:=c(x)\in(0,1)$ for which, since if $\left.\nabla_{V_{i}}G_{\epsilon}\right|_{x_{0}}:=A_{i}$,
$Q_{i}(A_{1},\ldots,A_{n})=0$, \begin{eqnarray*}
Q_{i} & = & \frac{\partial Q_{i}}{\partial P_{j}}(P_{1}(c),\ldots,P_{n}(c))(\nabla_{j}(G_{\epsilon}+t\eta)-A_{j})\\
 & = & \sum_{j\neq i}E\frac{\sqrt{\Pi_{k\neq i,j}\left(1+E\left\Vert P_{k}(c)\right\Vert ^{2}\right)}}{\sqrt{1+E\left\Vert P_{i}(c)\right\Vert ^{2}}\sqrt{1+E\left\Vert P_{j}(c)\right\Vert ^{2}}}<P_{j}(c),\nabla_{j}(G_{\epsilon}+t\eta)-A_{j}>\\
 &  & -E\frac{\sqrt{\Pi_{j\neq i}\left(1+E\left\Vert P_{j}(c)\right\Vert ^{2}\right)}}{\left(1+E\left\Vert P_{i}(c)\right\Vert ^{2}\right)^{3/2}}<P_{i}(c),\nabla_{i}(G_{\epsilon}+t\eta)-A_{i}>\end{eqnarray*}
for some $(P_{1}(c),\ldots,P_{n}(c))=(A_{1},\ldots,A_{n})+c(\nabla_{1}(G_{\epsilon}+t\eta)-A_{1},\ldots,\nabla_{n}(G_{\epsilon}+t\eta)-A_{n})$,
$c\in(0,1)$. 

Now, $f(t)=t/\sqrt{1+t}$ is increasing for $t>0$ and $\left\Vert P_{l}(c)\right\Vert \ll\left\Vert \nabla_{l}G_{\epsilon}\right\Vert \ll\left\Vert P_{l}(c)\right\Vert $
(which follows because $\left\Vert \nabla\eta\right\Vert $ and $A_{l}$
are bounded), so that \[
\sqrt{\Pi_{k\neq i,j}\left(1+E\left\Vert P_{k}(c)\right\Vert ^{2}\right)}\ll\sqrt{\Pi_{k\neq i,j}\left(1+E\left\Vert \nabla_{k}G_{\epsilon}\right\Vert ^{2}\right)}\]
 and \[
\frac{E<P_{j}(c),\nabla_{j}(G_{\epsilon}+t\eta)-A_{j}>}{\sqrt{1+E\left\Vert P_{j}(c)\right\Vert ^{2}}}\ll\frac{E\left\Vert P_{j}(c)\right\Vert ^{2}}{\sqrt{1+E\left\Vert P_{j}(c)\right\Vert ^{2}}}\ll\frac{E\left\Vert \nabla_{j}G_{\epsilon}\right\Vert ^{2}}{\sqrt{1+E\left\Vert \nabla_{j}G_{\epsilon}\right\Vert ^{2}}}.\]
 Then, applying these inequalities to the expression for $Q_{i}$
above,\begin{eqnarray*}
\left|Q_{i}\right| & \ll & \sum_{j}\frac{E\left\Vert \nabla_{j}G_{\epsilon}\right\Vert ^{2}\sqrt{\Pi_{k\neq i,j}\left(1+E\left\Vert \nabla_{k}G_{\epsilon}\right\Vert ^{2}\right)}}{\sqrt{1+E\left\Vert \nabla_{j}G_{\epsilon}\right\Vert ^{2}}\sqrt{1+E\left\Vert P_{i}(c)\right\Vert ^{2}}},\end{eqnarray*}
and so, this time because $f(t)=t/\sqrt{1+t^{2}}$ is also increasing,
and using Lemma (\ref{lem:lemma3c}) in the second step, \begin{eqnarray*}
\left|Q\right| & \ll & \frac{1}{\sqrt{E}}\int_{B(X_{0},1)\backslash\phi_{R}(D_{\epsilon})}\sum_{i,j}\frac{\sqrt{E}\left\Vert \nabla_{i}G_{\epsilon}\right\Vert E\left\Vert \nabla_{j}G_{\epsilon}\right\Vert ^{2}\sqrt{\Pi_{k\neq i,j}\left(1+E\left\Vert \nabla_{k}G_{\epsilon}\right\Vert ^{2}\right)}}{\sqrt{1+E\left\Vert \nabla_{i}G_{\epsilon}\right\Vert ^{2}}\sqrt{1+E\left\Vert \nabla_{j}G_{\epsilon}\right\Vert ^{2}}}d\mathcal{L}^{n}\\
 & = & \frac{1}{\sqrt{E}}\int_{B(X_{0},1)\backslash\phi_{R}(D_{\epsilon})}\sum_{i,j}\left.\frac{\left\Vert \nabla_{i}g_{\epsilon}\right\Vert \left\Vert \nabla_{j}g_{\epsilon}\right\Vert ^{2}\sqrt{\Pi_{k\neq i,j}\left(1+\left\Vert \nabla_{k}g_{\epsilon}\right\Vert ^{2}\right)}}{\sqrt{1+\left\Vert \nabla_{i}g_{\epsilon}\right\Vert ^{2}}\sqrt{1+\left\Vert \nabla_{j}g_{\epsilon}\right\Vert ^{2}}}\right|_{x=RX}d\mathcal{L}^{n}(X)\\
 & \ll & \frac{1}{\sqrt{E}}\int_{B(X_{0},1)\backslash\phi_{R}(D_{\epsilon})}\left(\left.\sqrt{\Pi_{k}\left(1+\left\Vert \nabla_{k}g_{\epsilon}\right\Vert ^{2}\right)}\right|_{x=RX}-1\right)d\mathcal{L}^{n}(X)\\
 & = & \sqrt{E}.\end{eqnarray*}
The last inequality follows from the fact that \[
4(\sqrt{1+a^{2}}\sqrt{1+b^{2}}c-1)-\frac{ab^{2}c}{\sqrt{1+a^{2}}\sqrt{1+b^{2}}}\geq0\]
for any $c>1$, which is a straightforward calculation.

On the bad set $D_{\epsilon}$, by the strong approximation theorem
\cite[4.2.20]{GMT} we can assume without loss of generality that
$T_{\epsilon}\rest\pi^{-1}(D_{\epsilon})$ is the image $\psi_{\#}(P)$,
where $P$ is a polyhedral chain and $\psi$ is Lipschitz. The definition
of $K(t)$ and the fact that the deformation $h_{t}$ is vertical
{[}cf. (\ref{sub:Euler-Lagrange-equations-for-T-epsilon}){]} implies
that \begin{eqnarray*}
\frac{d}{dt}K(t) & = & \left.\frac{d}{dt}\right|_{t}\int_{\pi^{-1}\left(D_{\epsilon}\right)}f_{\epsilon}(\overrightarrow{T_{t}})d\left\Vert T_{t}\right\Vert /(ER^{n})\\
 & = & \int_{\pi^{-1}(D_{\epsilon})}\frac{d}{dt}d\left\Vert T_{t}\right\Vert /(ER^{n})\end{eqnarray*}
since the deformation will leave the penalty part fixed. In addition,
the derivative of this integrand will be 0 at all points with a vertical
tangent plane, again due to the fact that the deformation is vertical.
At all points where the tangent plane is not vertical, the mean-value
theorem approximation used for the good set will again hold, where
we can replace $g_{\epsilon}(x)$ by $\psi(p)$, where $\pi(\psi(p))=x.$
In the notation above, if \[
\mathcal{F}_{\epsilon,R}(S):=\left\Vert (\phi_{R}^{-1})_{\#}(S)\right\Vert /(ER^{n})+\frac{1}{\epsilon ER^{n}}H(S),\]
then applying the squash-deformation, for which $\phi_{R}\psi:=\Psi$\begin{eqnarray*}
\frac{d}{dt}K(t) & = & \frac{d}{dt}\mathcal{F}_{\epsilon,R}\left((H_{t})_{\#}(\Psi_{\# }(P))\rest\left(\phi_{R}(\pi^{-1}(D_{\epsilon}))\right)\right)\\
 & = & \frac{1}{ER^{n}}\frac{d}{dt}\int_{P}\sqrt{\sum_{|\alpha|+|\beta|=n}E^{|\beta|}((H_{t})_{\#}(\Psi_{\# }(P))_{\alpha\beta})^{2}}d\left\Vert P\right\Vert ,\end{eqnarray*}
where again the penalty part is irrelevant since the deformation is
vertical, and the deformation $H_{t}$ defined by $H_{t}=\phi_{R}h_{t}\phi_{R}^{-1}$
becomes translation vertically by $t\eta(X)$, where $X=\pi(p)$,
$p\in Supp(P)$. Also as a consequence of the verticality of the deformation,
the $\beta=0$ term of the integral will be unchanged under the deformation,
so \begin{eqnarray*}
\frac{d}{dt}K(t) & = & \frac{d}{dt}\mathcal{F}_{\epsilon,R}\left((H_{t})_{\#}(\Psi_{\# }(P))\rest\left(\phi_{R}(\pi^{-1}(D_{\epsilon}))\right)\right)\\
 & = & \frac{1}{ER^{n}}\int_{P}\frac{E\sum_{|\alpha|+|\beta|=n,\beta\neq0}E^{|\beta|-1}(H_{t})_{\#}(\Psi_{\# }(P))_{\alpha\beta}\frac{d}{dt}(H_{t})_{\#}(\Psi_{\# }(P))_{\alpha\beta}}{\sqrt{\sum_{|\alpha|+|\beta|=n}E^{|\beta|}((H_{t})_{\#}(\Psi_{\# }(P))_{\alpha\beta})^{2}}}d\left\Vert P\right\Vert (p)\\
 & = & \frac{1}{R^{n}}\left(\int_{D_{\epsilon}}\sum_{i}\left(<\nabla_{i}F,\nabla_{i}\eta>+t<\nabla_{i}\eta,\nabla_{i}\eta>\right)\mathcal{A}_{ii}d\mathcal{L}^{n}+Q\right).\end{eqnarray*}
The factors $Q_{i}$, $Q=\int_{D_{\epsilon}}\sum_{i}\left(<\nabla_{i}\eta,\nabla_{i}F+t\nabla_{i}\eta>\right)Q_{i}d\mathcal{L}^{n}$
can be bounded as before. The factorization of the integrand \[
\sqrt{\sum_{|\alpha|+|\beta|=n}E^{|\beta|}((H_{t})_{\#}(\Psi_{\# }(P))_{\alpha\beta})^{2}}=\sqrt{\prod_{j}\left(1+E\left\Vert P_{j}\right\Vert ^{2}\right)},\]
since we only are concerned with points at non-vertical tangents,
$P_{j}=\nabla_{j}(F+t\eta)$, is well-defined, where the covariant
derivative is in the direction of $\partial/\partial X_{j}$ as before,
and the basis is chosen to diagonalize the quadratic form $(V,W)\mapsto<\nabla_{V}F+t\eta,\nabla_{W}F+t\eta>$
as in the previous case, $A$ is this quadratic form at $t=0$, and
$\mathcal{A}$ is derived from $A$ as before. Each such $Q_{i}$
can also be bounded as (since $E<1$) by\begin{eqnarray*}
\left|\left(\left\langle \nabla_{i}\eta,\nabla_{i}F+t\nabla_{i}\eta\right\rangle \right)Q_{i}\right| & \ll & \sqrt{E}\sqrt{\Pi_{j\neq i}\left(1+E\left\Vert P_{j}\right\Vert ^{2}\right)}\\
 & \ll & \sqrt{E}\sqrt{\prod_{j}\left(1+E\left\Vert P_{j}\right\Vert ^{2}\right)}\\
 & = & \sqrt{E}\sqrt{\sum_{|\alpha|+|\beta|=n}E^{|\beta|}((H_{t})_{\#}(\Psi_{\# }(P))_{\alpha\beta})^{2}},\end{eqnarray*}
 so that \begin{eqnarray*}
 &  & \left|\frac{d}{dt}K(t)-\int_{D_{\epsilon}}\sum_{i,k}\mathcal{A}_{ik}\left\langle \nabla_{i}\eta,\nabla_{k}F+t\nabla_{k}\eta\right\rangle d\mathcal{L}^{n}\right|\\
 & \ll & \frac{1}{ER^{n}}\int_{P}\frac{d}{dt}\sqrt{\sum_{|\alpha|+|\beta|=n}E^{|\beta|}((H_{t})_{\#}(\Psi_{\# }(P))_{\alpha\beta})^{2}}d\left\Vert P\right\Vert \\
 & \ll & \left\Vert T_{\epsilon}\rest\pi^{-1}(D_{\epsilon})\right\Vert /(\sqrt{E}R^{n})\\
 & \ll & \sqrt{E}\end{eqnarray*}
by Corollary (\ref{lem:lemma 3a-improved}). This establishes the
Lemma.
\hfill{}$\square$\end{pf}
\begin{lem}
\label{lem:lemma 13} With the hypotheses of Lemma (\ref{lemma 12}),
if the support of $\eta$ is contained in $|X|<1-E^{1/4n}$ and $|\nabla\eta|\leq1$,
we also have\[
\left|\int_{B(X_{0},1)}\sum\mathcal{A}_{ik}\left\langle \nabla_{i}\eta,\nabla_{k}F\right\rangle d\mathcal{L}^{n}\right|\ll\sqrt{E}.\]

\end{lem}
\begin{pf}
Here we use the minimality of $T_{\epsilon}$. From Lemma (\ref{lemma 12}),
we have that \[
\left|\frac{d}{dt}\mathcal{F}_{\epsilon,R}\left(\left(h_{t}\right)_{\#}\left(T_{\epsilon,R}\right)\right)/(ER^{n})-\int_{B(X_{0},1)}\sum_{i,k}\mathcal{A}_{ik}\left\langle \nabla_{i}\eta,\nabla_{k}F+t\nabla_{k}\eta\right\rangle d\mathcal{L}^{n}\right|\ll\sqrt{E}.\]
However, since $T_{\epsilon}$ minimizes $\mathcal{F}_{\epsilon}$,
$T_{\epsilon,R}$will minimize $\mathcal{F}_{\epsilon,R}$ by its
definition. This implies that \[
\left.\frac{d}{dt}\right|_{0}\mathcal{F}_{\epsilon,R}\left(\left(h_{t}\right)_{\#}\left(T_{\epsilon,R}\right)\right)=0,\]
and the Lemma follows from setting $t=0$.
\hfill{}$\square$\end{pf}
\begin{lem}
\label{lem:lemma 14}For any $L:B(x_{0},R)\rightarrow\R^{k}$ so that,
for some $\sigma$, $\left|grad(L)\right|\leq\sigma\leq1$. , let
$h(x,y)=(x,y-L(x))$. Then \[
Exc(h_{\#}(T),R)\ll E+\sigma^{2}.\]

\end{lem}
\begin{pf}
Since $h$ is vertical, if $\mathbf{e}=dx^{1}\wedge\cdots\wedge dx^{n}$
is the horizontal $n$-vector in $\Lambda_{n}(B(x_{0},R)\times\R^{k})$,
\[
<\mathbf{e},h_{\#}\left(\overrightarrow{T_{\epsilon}}\right)>=<\mathbf{e},\overrightarrow{T_{\epsilon}}>\]
and so, for any multiindex\[
\left|<dx^{\alpha}\wedge dy^{\beta},h_{\#}\left(\overrightarrow{T_{\epsilon}}\right)>-<dx^{\alpha}\wedge dy^{\beta},\overrightarrow{T_{\epsilon}}>\right|\ll\sigma.\]
Since \begin{eqnarray*}
\left\Vert h_{\#}(\overrightarrow{T})\right\Vert  & = & \sqrt{\sum_{|\alpha|+|\beta|=n}<dx^{\alpha}\wedge dy^{\beta},h_{\#}(\overrightarrow{T})>^{2}}\\
 & \leq & \sqrt{<\mathbf{e},\overrightarrow{T}>+\sum_{|\alpha|+|\beta|=n,|\beta|>0}(<dx^{\alpha}\wedge dy^{\beta},\overrightarrow{T}>+c\sigma)^{2}}\\
 & \leq & \sqrt{1+c'\sigma\left(\sum_{|\alpha|+|\beta|=n,|\beta|>0}<dx^{\alpha}\wedge dy^{\beta},\overrightarrow{T}>\right)+c'\sigma^{2}}\\
 & \leq & \sqrt{1+c''\sigma\sqrt{\sum_{|\alpha|+|\beta|=n,|\beta|>0}(<dx^{\alpha}\wedge dy^{\beta},\overrightarrow{T}>)^{2}}+c''\sigma^{2}}\\
 & \leq & 1+c''\sigma\sqrt{1-<\mathbf{e},\overrightarrow{T}>^{2}}+c''\sigma^{2},\end{eqnarray*}
\begin{eqnarray*}
\left\Vert h_{\#}(T)\right\Vert  & \leq & (1+c''\sigma^{2})\left\Vert T\right\Vert +c''\sigma\int_{C(x_{0},R)}\sqrt{1-<\mathbf{e},\overrightarrow{T}>^{2}}d\left\Vert T\right\Vert \\
 & \leq & (1+c''\sigma^{2})\left\Vert T\right\Vert +c''\left(\sigma\sqrt{\left\Vert T\right\Vert }\right)\sqrt{ER^{n}}\\
 & \leq & \left\Vert T\right\Vert +c'''\sigma^{2}\left\Vert T\right\Vert +c'''ER^{n}.\end{eqnarray*}
Since $\left\Vert T\right\Vert \ll R^{n}$, the Lemma follows.
\hfill{}$\square$\end{pf}

\section{Iterative inequality}

Fix $\beta$, $0<\beta\leq1/4$. 

\begin{prop}
\label{lem:Main Lemma}If $T$ is a mass-minimizing rectifiable section
$T\in\widetilde{\Gamma}(B)$ which is the limit of a sequence of penalty
minimizers $T_{\epsilon},$ and there exists a positive constant $\alpha=\alpha(\beta)$
and a constant $c$, so that if \[
R+Exc(T;R)\leq\alpha,\]
then \begin{equation}
Exc(h_{\#}T;\beta R)\leq c\beta^{2}Exc(T;R)\label{eq:main-lemma-1}\end{equation}
for some linear map $h(x,y)=(x,y-l(x))$ with \begin{equation}
|\textrm{grad $l|\leq\alpha^{-1}\sqrt{Exc(T;R)}.$}\label{eq:main-lemma-2}\end{equation}

\end{prop}
\begin{rem}
Note that, if this Lemma holds with some one value of $\alpha$, it
will also hold with any smaller $\alpha$. Also, recall from Theorem
{[}\ref{thm:penalty argument}{]} that for any homology class of rectifiable
sections there will be one such section which is the limit of penalty
minimizers.
\end{rem}
\begin{pf}
If this is not the case, then we will be able to find a sequence $R_{i}\rightarrow0$,
$\epsilon_{i}\rightarrow0$, along with functionals $\mathcal{F}_{i}:=\mathcal{F}_{\epsilon_{i},R_{i}}$
as above and $T_{i}\rightarrow T$ (minimizers of $\mathcal{F}_{i}$),
and excesses $E_{i}:=Exc(T_{\epsilon_{i}};R_{i},x_{0})$ for which
$E_{i}\rightarrow0$ and (by choosing each $R_{i}$ sufficiently small)
$E_{i}^{1/4n}/\epsilon_{i}\rightarrow0$, and \begin{equation}
\limsup_{i\rightarrow\infty}E_{i}^{-1}Exc((h_{i})_{\#}(T_{i});\beta R_{i})\geq c\beta^{2}\label{eq:non-main-lemma-1}\end{equation}
for all linear maps $h_{i}(x,y)=(x,y-l_{i}(x))$with \[
\limsup_{i\rightarrow\infty}E_{i}^{-1/2}|\textrm{grad $l_{i}|<\infty.$}\]
Such a sequence $\{ T_{i},\,\mathcal{F}_{i},\, R_{i}\}$, following
\cite{Bombieri}, will be called an \emph{admissible sequence}.

As before, let $D_{\epsilon_{i}}$ be the bad set over which $T_{i}:=T_{\epsilon_{i}}$
is not necessarily a $C^{1}$ graph with bounded gradient, and let
$D_{i}:=\phi_{R_{i}}(D_{\epsilon_{i}})\cap B(X_{0},1)$. Then, on
$B(X_{0},1)\backslash D_{i}$, $T_{i}:=T_{\epsilon_{i},R_{i}}$will
be the graph of a $C^{1}$ function $G_{i}$, agreeing on $B(X_{0},1)\backslash D_{i}$
with $F_{i}$, which is the BV carrier of $T_{i}$ on $B(X_{0},R)$.
We need to show:
\begin{lem}
\label{pro:lemma 15}For all $i$ sufficiently large
\begin{enumerate}
\item \[
\int_{B(X_{0},1)}\left\Vert dF_{i}\right\Vert d\mathcal{L}^{n}\ll1,\]

\item \[
\lim_{i}\left\Vert D_{i}\right\Vert =0\]

\item \[
\lim_{i}\frac{\int_{B(X_{0},1/2)\backslash\phi_{R_{i} }(D_{i})}\left|F_{i}\right|^{2p}d\mathcal{L}^{n}}{\left(E_{i}\right)^{p/2n}}\ll_{p}1,\,\,1\leq p<\frac{n}{n-1},\]

\item \[
\int_{B(X_{0},1)}\left|F_{i}\right|d\mathcal{L}^{n}\ll1\]

\item \[
\frac{Exc(T_{i},R_{i}/4)}{E_{i}}\ll\left(2+\frac{1}{2\epsilon_{i}}\right)E_{i}^{1/4n}+\left(2+\frac{1}{2\epsilon_{i}}\right)^{3/2}E_{i}^{3/4n},\]

\item The limit\[
\lim_{i}\mathcal{A}_{i}:=\mathcal{A}_{0}\]
 is the symbol of an elliptic PDE. 
\item for every smooth $\eta(X)$ with compact support in $|X|<1$ we have
\[
\lim_{i}\int_{B(x_{0},1)}\sum\left(A_{i}\right)_{jk}\left\langle \frac{\partial\eta}{\partial X_{j}},D_{k}F_{i}\right\rangle d\mathcal{L}^{n}=0,\]

\item Finally, if $h_{i}(x,y)=(x,y-l_{i}(x))$ is a sequence of linear maps
with \[
\lim_{i}\frac{|grad(l_{i})|}{\sqrt{E_{i}}}\leq\sigma\]
then \[
\lim_{i}\frac{Exc((h_{i})_{\#}(T_{i}),R_{i})}{E_{i}}\ll(1+\sigma^{2}).\]

\end{enumerate}
\end{lem}
\begin{pf}
Set, for each $i$ in the sequence, \[
\overline{y}(i)\,:=\frac{1}{\left\Vert B(x_{0},R/2)\backslash\left(D_{\epsilon}\cap B(x_{0},R/2)\right)\right\Vert }\int_{B(x_{0},R/2)\backslash\left(D_{\epsilon}\cap B(x_{0},R/2)\right)}f_{i}d\mathcal{L}^{n}.\]
For each $i$, translate the corresponding graph so that so that $\overline{y}(i)=\overline{0}.$
By Lemma(\ref{lem:lemma7}), there is a constant $\tau$ so that for
all $p$, $1\leq p\leq\frac{n}{n-1}$,\[
\left(\frac{\int_{B(x_{0},R)}\left|f_{i}\right|^{p}d\mathcal{L}^{n}}{R^{n}}\right)^{1/p}\leq\tau\frac{\int_{B(x_{0},R)}\left\Vert df_{i}\right\Vert d\mathcal{L}^{n}}{R^{n-1}}:=\tau\frac{\int_{B(x_{0},R)}\left\Vert df_{i}\right\Vert }{R^{n-1}},\]
and since we have by Lemma (\ref{lem:lemma2}) that $\int_{B(x_{0},R)}\left\Vert df_{i}\right\Vert \ll\sqrt{E_{i}}R^{n}$,
with $p=1$we conclude that \[
\int_{B(x_{0},R)}\left|f_{i}\right|d\mathcal{L}^{n}\ll\sqrt{E_{i}}R^{n+1},\]
which since $F_{i}(X)=f_{i}(RX)/(\sqrt{E_{i}}R),$ as before yields
that for all $i$ sufficiently large, \[
\int_{B(X_{0},1)}\left\Vert dF_{i}\right\Vert d\mathcal{L}^{n}\ll1,\textrm{ and $\int_{B(X_{0},1)}\left|F_{i}\right|d\mathcal{L}^{n}\ll1$,}\]
 which are statements (4) and (1), respectively. Lemma (\ref{lem:lemma 8})
and the definition of $F_{i}$ immediately gives statement (3), and
statement (2) follows from the bound $\left\Vert D_{\epsilon}\right\Vert \leq ER^{n}$,
so that $\left\Vert D_{i}\right\Vert \ll E_{i}$, and the choice of
$E_{i}\downarrow0$.

To show statement (5), use Lemma(\ref{lemma_11}) to show that (with
$\overline{y}=0$) and Lemma (\ref{lem:lemma 8})

\begin{eqnarray*}
\frac{Exc(T_{i},R_{i}/4)}{E_{i}} & \ll & (1+\frac{1}{2\epsilon_{i}})\left(\mu+E_{i}^{1/2n}\left(\frac{1}{\mu}\right)\right)+\left.\int_{B(x_{0},R_{i}/2)\backslash D_{\epsilon_{i}}}\left|f_{i}\right|^{2}d\mathcal{L}^{n}\right/(\mu E_{i}R_{i}^{n+2})\\
 & \ll & (1+\frac{1}{2\epsilon_{i}})\left(\mu+E_{i}^{1/2n}\left(\frac{1}{\mu}\right)\right)+\left.\int_{B(x_{0}1/2)}\left|F_{i}\right|^{2}d\mathcal{L}^{n}\right/(\mu)\\
 & \ll & (1+\frac{1}{2\epsilon_{i}})\left(\mu+E_{i}^{1/2n}\left(\frac{1}{\mu}\right)\right)+\left.\int_{B(x_{0}1/2)}\left|F_{i}\right|^{2}d\mathcal{L}^{n}\right/(\mu)\\
 & = & \mu(1+\frac{1}{2\epsilon_{i}})+\frac{1}{\mu}\left((1+\frac{1}{2\epsilon_{i}})E_{i}^{1/2n}+\int_{B(x_{0}1/2)}\left|F_{i}\right|^{2}d\mathcal{L}^{n}\right).\end{eqnarray*}
Taking $\mu$ to minimize the right hand side above, \[
\mu=\mu_{i}=\frac{\sqrt{(1+\frac{1}{2\epsilon_{i}})E_{i}^{1/2n}+\int_{B(x_{0}1/2)}\left|F_{i}\right|^{2}d\mathcal{L}^{n}}}{\sqrt{1+\frac{1}{2\epsilon_{i}}}},\]
 which for $i$ sufficiently large will be less than one, by (3) above,
and the fact that $E_{i}\searrow0$, gives \begin{eqnarray*}
\frac{Exc(T_{i},R_{i}/4)}{E_{i}} & \ll & \sqrt{1+\frac{1}{2\epsilon_{i}}}\sqrt{(1+\frac{1}{2\epsilon_{i}})E_{i}^{1/2n}+\int_{B(x_{0}1/2)}\left|F_{i}\right|^{2}d\mathcal{L}^{n}}\\
 &  & +\left((1+\frac{1}{2\epsilon_{i}})E_{i}^{1/2n}+\int_{B(x_{0}1/2)}\left|F_{i}\right|^{2}d\mathcal{L}^{n}\right)^{3/2}\\
 & \ll & \sqrt{1+\frac{1}{2\epsilon_{i}}}\sqrt{2+\frac{1}{2\epsilon_{i}}}E^{1/4n}+\left(2+\frac{1}{2\epsilon}\right)^{3/2}E^{3/4n}.,\end{eqnarray*}
easily giving (5).

Statement (6) follows from Equation (\ref{eq:A_ij}). Statement (7)
follows from Lemma (\ref{lem:lemma 13}).

Finally, statement (8) follows from Lemma (\ref{lem:lemma 14}).
\hfill{}$\square$\end{pf}
By statements (1) and (3) of this Lemma, invoking the closure and
compactness theorems for BV functions \cite{GMT}, we can assume that
there is an element $F\in BV(B(X_{0},1))$ so that a subsequence (which
by standard abuse of notation we do not re-label) $F_{i}\rightarrow F$
strongly in $L^{1}(B(X_{0},1))$ and $DF_{i}\rightarrow DF$ as distributions.
We then have \[
\int_{B(X_{0},1)}\sum A_{jk}\left\langle \frac{\partial\eta}{\partial X_{j}},D_{k}F\right\rangle d\mathcal{L}^{n}=0,\]
for all smooth $\eta$ with compact support in $\left|X\right|<1$.
Thus, $F$ will be $\mathcal{A}$-harmonic, and thus is a real-analytic
function. It then follows from the Di Giorgi-Moser-Morrey estimates
for diagonal elliptic systems \cite{Morrey} that\[
\sup_{B(X_{0},1/2)}\left|F\right|\ll\int_{B(X_{0},1)}\left|F\right|d\mathcal{L}^{n}=\lim_{i}\int_{B(X_{0},1)}\left|F_{i}\right|d\mathcal{L}^{n}\ll1,\]
so we can shift the graph so that $F(X_{0})=0$. Our previous shift
was chosen so that, for each $i$, $\overline{y}(i)=0$, where \[
\overline{y}(i):=\frac{1}{\left\Vert B(x_{0},R_{i}/2)\backslash\left(D_{\epsilon_{i}}\cap B(x_{0},R_{i}/2)\right)\right\Vert }\int_{B(x_{0},R_{i}/2)\backslash\left(D_{\epsilon_{i}}\cap B(x_{0},R_{i}/2)\right)}f_{i}d\mathcal{L}^{n}.\]
The bounds on the $L^{1}$ norms of $F$ and the BV-norm of $DF$
are not worsened by this assumption except for possible change of
constants, which are implicit in the notation. In addition, the bound
of statement (3) in Lemma(\ref{pro:lemma 15}) continues to hold as
well, since the bound $\overline{y}(i)=0$ becomes \[
\int_{B(X_{0},1/2)\backslash\phi_{R_{i} }(D_{i})}F_{i}d\mathcal{L}^{n}=0,\]
from which follows the fact that \[
\int_{B(X_{0},1/2)\backslash\phi_{R_{i} }(D_{i})}\left|F_{i}\right|^{2}d\mathcal{L}^{n}\leq\int_{B(X_{0},1/2)\backslash\phi_{R_{i} }(D_{i})}\left|F_{i}+C\right|^{2}d\mathcal{L}^{n}\]
for any constant vector $C$. 

\begin{lem}
\label{lem:Lemma 16} Let $\{ T_{i},\,\mathcal{F}_{i},\, R_{i}\}$
be admissible. Under a suitable translation (or change of coordinates),
$F(0)=0,\, F_{i}\rightarrow F$ strongly in $L^{1}$, $F$ is a solution
to the equation \[
\int_{B(X_{0},1)}\sum\left(\mathcal{A}\right)_{jk}\left\langle \frac{\partial\eta}{\partial X_{j}},D_{k}F\right\rangle d\mathcal{L}^{n}=0,\]
as well as \[
\int_{B(X_{0},1)}\left|F\right|d\mathcal{L}^{n}+\int_{B(X_{0},1)}\left|gradF\right|d\mathcal{L}^{n}\ll1,\]
\[
\sup_{B(X_{0},1/2)}\left(\left|F\right|,\left|gradF\right|\right)\ll1,\]
 and\[
\lim_{i}\frac{Exc(T_{i},R_{i}/4)}{E_{i}}=0\]
\textcolor{black}{under the assumption that $\lim_{i}E_{i}^{1/4n}/\epsilon_{i}=0$.}
\end{lem}
\begin{pf}
\[
\int_{B(X_{0},1)}\left|F\right|d\mathcal{L}^{n}=\lim_{i}\int_{B(X_{0},1)}\left|F_{i}\right|d\mathcal{L}^{n}\]
because $F_{i}\rightarrow F$ strongly in $L^{1}$, and \[
\int\left\Vert DF\right\Vert d\mathcal{L}^{n}\leq\lim_{i}\int_{B(X_{0},1)}\left\Vert DF_{i}\right\Vert d\mathcal{L}^{n}\]
 by lower semi-continuity with respect to BV-convergence. In order
to complete the proof of the Lemma, we need only show that \[
\lim_{i}\int_{B(X_{0},1/2)\backslash\phi_{R_{i}}(D_{i})}\left|F_{i}\right|^{2}d\mathcal{L}^{n}=\int_{B(X_{0},1/2)}\left|F\right|^{2}d\mathcal{L}^{n}.\]
But, if $\Phi_{i}$ is the characteristic function of $B(X_{0},1)\backslash\phi_{R_{i}}(D_{i})$,
then \begin{eqnarray*}
\lim_{i}\left|\int_{B(X_{0},1/2)\backslash\phi_{R_{i}}(D_{i})}\left|F_{i}\right|^{2}d\mathcal{L}^{n}\right.\\
\left.-\int_{B(X_{0},1/2)}\left|F\right|^{2}d\mathcal{L}^{n}\right| & = & \lim_{i}\left|\int_{B(X_{0},1/2)}\Phi_{i}\left|F_{i}\right|^{2}-\left|F\right|^{2}d\mathcal{L}^{n}\right|\\
 & = & \lim_{i}\left|\int_{B(X_{0},1/2)}\Phi_{i}\left(\left|F_{i}\right|^{2}-\left|F\right|^{2}\right)+\Phi_{i}\left|F\right|^{2}-\left|F\right|^{2}d\mathcal{L}^{n}\right|\\
 & = & \lim_{i}\left|\int_{B(X_{0},1/2)}\Phi_{i}\left(\left|F_{i}\right|-\left|F\right|\right)\left(\left|F_{i}\right|+\left|F\right|\right)+\left(\Phi_{i}-1\right)\left|F\right|^{2}d\mathcal{L}^{n}\right|\\
 & \leq & \lim_{i}\left|\int_{B(X_{0},1/2)}\Phi_{i}\left|F_{i}-F\right|\left(\left|F_{i}\right|+\left|F_{i}\right|\right)+\left(\Phi_{i}-1\right)\left|F\right|^{2}d\mathcal{L}^{n}\right|\\
 & \leq & \lim_{i}\int_{B(X_{0},1/2)}\Phi_{i}\left(\left|F_{i}\right|+\left|F\right|\right)\left|F_{i}-F\right|d\mathcal{L}^{
}\\
 &  & +\int_{B(X_{0},1/2)}(1-\Phi_{i})\left|F\right|^{2}d\mathcal{L}^{n}.\end{eqnarray*}
Since $F$ is uniformly bounded in $B(X_{0},1/2)$ and $F_{i}\rightarrow F$
strongly in $L^{1}$, a subsequence will converge almost-everywhere
pointwise, and $\lim_{i}\int_{B(X_{0},1)}(1-\Phi_{i})d\mathcal{L}^{n}=0$,
the last integral above goes to 0, and \begin{eqnarray*}
\lim_{i}\left|\int_{B(X_{0},1/2)}\Phi_{i}\left|F_{i}\right|^{2}-\left|F\right|^{2}d\mathcal{L}^{n}\right| & \leq & 2\lim_{i}\int_{B(X_{0},1/2)}\Phi_{i}\left|F_{i}\right|\left|F_{i}-F\right|d\mathcal{L}^{
}\\
 & \leq & \lim_{i}\left(\int_{B(X_{0},1/2)}\Phi_{i}\left|F_{i}\right|^{2p-1}\left|F_{i}-F\right|d\mathcal{L}^{
}\right)^{\frac{1}{2p-1}}\cdot\\
 &  & \cdot\left(\int_{B(X_{0},1/2)}\left|F_{i}-F\right|d\mathcal{L}^{
}\right)^{1-\frac{1}{2p-1}}.\end{eqnarray*}
The last step is H\"{o}lder's inequality for the measure $\mu=|F_{i}-F|d\mathcal{L}^{n}$.
If $1<p<\frac{n}{n-1}$ then the first of these last two integrals
is uniformly bounded by statement (3) of Lemma (\ref{pro:lemma 15})
by a power of $E_{i}$ (Note that $E_{i}\rightarrow0$ as $i\rightarrow\infty$.),
and the height bound on $F_{i}$ and $F$ coming from the compactness
of the fiber of the bundle. The second integral goes to 0 in the limit
by the strong convergence of $F_{i}$ to $F$ in $L^{1}$.
\hfill{}$\square$\end{pf}
We can now complete the proof of Proposition (\ref{lem:Main Lemma}).
Let $L=L(X)$ denote the linear forms \[
L(X):=\sum\frac{\partial F}{\partial X_{i}}(0)X_{i},\]
and let $h_{i}$ be the maps \[
h_{i}(x,y)=\left(x,y-\sqrt{E_{i}}L(x)\right),\]
and \[
\widetilde{T_{i}}:=(h_{i})_{\#}(T_{i}).\]
Set $\widetilde{E_{i}}:=Exc(\widetilde{T_{i}},R_{i})$. Since $|gradL|=\left|DF\right|\ll1$,
we apply statement (8) of Lemma (\ref{pro:lemma 15}), which shows
that \[
\lim_{i}\frac{\widetilde{E_{i}}}{E_{i}}\ll1.\]
\setcounter{case}{0}

\begin{case}
$\lim_{i}\widetilde{E_{i}}/E_{i}=0.$ This contradicts \[
\limsup_{i\rightarrow\infty}E_{i}^{-1}Exc((h_{i})_{\#}(T_{i});\beta R_{i})\geq c\beta^{2},\]
which is a basic assumption on the sequence $T_{i}$, since this case
implies that $\lim_{i}Exc(\widetilde{T_{i}},\beta R_{i})/E_{i}\leq\lim_{i}\widetilde{E_{i}}/(\beta^{n}E_{i})=0.$
\end{case}

\begin{case}
$\lim_{i}\widetilde{E_{i}}/E_{i}>0.$ Then, the currents $\widetilde{T_{i}}$
minimize $\mathcal{F}_{i}$, so that $\{\widetilde{T_{i}},\mathcal{F}_{i},R_{i}\}$
will still be admissible, so that \[
\widetilde{F_{i}}\rightarrow\widetilde{F},\]
where \[
\widetilde{F_{i}}(X):=\sqrt{\widetilde{E_{i}}^{-1}E_{i}}\left(F_{i}-L\right)(X).,\,\,\widetilde{F}(X):=\lim_{i}\sqrt{\widetilde{E_{i}}^{-1}E_{i}}\left(F-L\right)(X).\]
Since $\widetilde{F}$ satisfies the conditions of Lemma (\ref{lem:Lemma 16}),
in particular \[
\lim_{i}\widetilde{E_{i}}^{-1}Exc(\widetilde{T_{i}},R_{i}/4)=0\]
 and in addition we have \[
D\widetilde{F}(0)=0,\]
 and the inequality (\ref{eq:non-main-lemma-1}) in the beginning
of the proof Proposition (\ref{lem:Main Lemma}), becomes \[
\lim_{i}E_{i}^{-1}Exc(\widetilde{T_{i}},\beta R_{i})\geq c\beta^{2}.\]
Define $s$ as the integer so that \[
\frac{1}{4}\leq4^{s}\beta<1\]
(assume that $\beta<1/4$, so that $s\geq1$), and for $\sigma=0,1,2,\ldots,s$
we consider \[
\widetilde{E_{i}}^{(\sigma)}:=Exc(\widetilde{T_{i}},4^{\sigma}\beta R_{i}).\]
It is clear by the fact that $\lim_{i}\frac{\widetilde{E_{i}}}{E_{i}}\ll1$
that\[
\widetilde{E_{i}}^{(\sigma)}\leq(4^{\sigma}\beta)^{-n}\widetilde{E_{i}}\ll(4^{\sigma}\beta)^{-n}E_{i},\]
that is, for some constant $C$, \[
\widetilde{E_{i}}^{(\sigma)}\leq C(4^{\sigma}\beta)^{-n}E_{i}.\]
If, for some $\sigma$ we have \[
\lim_{i}E_{i}^{-1}\widetilde{E_{i}}^{(\sigma)}=0,\]
then we have \[
\lim_{i}E_{i}^{-1}Exc(\widetilde{T_{i}},\beta R_{i})=0,\]
contradicting our assumption above. So, we can assume that \[
\lim_{i}E_{i}^{-1}\widetilde{E_{i}}^{(\sigma)}>0.\]
We also have \[
\lim_{i}E_{i}^{-1}\widetilde{E_{i}}^{(\sigma)}\ll(4^{\sigma}\beta)^{-n}<+\infty,\]
for $\sigma=0,1,\ldots,s$. Now, these inequalities and the fact that
$\{\widetilde{T_{i}},\mathcal{F}_{i},R_{i}\}$ is admissible implies
that also \[
\{\widetilde{T_{i}}\rest C(x_{0},4^{\sigma}\beta R_{i}),\mathcal{F}_{i}\circ h_{\#}^{-1},4^{\sigma}\beta R_{i}\}\]
will be admissible. Thus, by the conclusions of Lemma (\ref{lem:Lemma 16}),\[
\lim_{i}\frac{\widetilde{E_{i}}^{(\sigma-1)}}{\widetilde{E_{i}}^{(\sigma)}}=0,\]
or, given any $a>0$, for $i$ sufficiently large, \[
\frac{\widetilde{E_{i}}^{(\sigma-1)}}{\widetilde{E_{i}}^{(\sigma)}}<a.\]
Iterating this inequality,\[
Exc(\widetilde{T_{i}},\beta R_{i}):=\widetilde{E_{i}}^{(0)}<a\widetilde{E_{i}}^{(1)}<\cdots<a^{\sigma}\widetilde{E_{i}}^{(\sigma)}\leq Ca^{\sigma}(4^{\sigma}\beta)^{-n}E_{i}.\]
Choosing $a$ sufficiently small will then guarantee that, for $i$
sufficiently large \[
\frac{Exc(\widetilde{T_{i}},\beta R_{i})}{E_{i}}<c\beta^{2},\]
contradicting the assumption, and completing the proof of Proposition
(\ref{lem:Main Lemma}). 
\end{case}
\hfill{}$\square$\end{pf}
There is a small extension of this Proposition that will be needed
for its application:

\begin{cor}
\label{cor:Main-Lemma-tilted}Given $\beta,$ $T$, $\alpha(\beta)$
as in Proposition (\ref{lem:Main Lemma}), then the conclusion of
the Proposition will still hold, for some $\alpha>0$, for the current
$H_{\#}(T)$, where $H(x,y)=(x,y+L(x))$ is a fixed linear map. That
is, if \begin{equation}
R+Exc(H_{\#}(T),R)\leq\alpha,\label{eq:excess bound}\end{equation}
then \begin{equation}
Exc(h_{\#}(H_{\#}(T)),\beta R)\leq c\beta^{2}Exc(H_{\#}(T),R)\label{eq:iterative excess bound}\end{equation}
for some linear map $h(x,y)=(x,y-l(x))$ with \begin{equation}
|\textrm{grad $l|\leq\alpha^{-1}\sqrt{Exc(H_{\#}(T),R)}.$}\label{eq:gradient bound}\end{equation}

\end{cor}
\begin{pf}
This corollary will follow from Proposition (\ref{lem:Main Lemma})
once it is shown that, under these conditions, $|\textrm{grad }l|$
is bounded as indicated in (\ref{eq:gradient bound}). However, under
the assumptions on $R$ and the excess of $H_{\#}(T)$, if $f$ is
the BV-carrier of $T$, so that $(f+L)$ is the BV-carrier of $H_{\#}(T)$,\begin{eqnarray*}
\int_{B(0,R)}\left\Vert D(f+L)\right\Vert d\mathcal{L}^{n} & \geq & \int\left|\left\Vert Df\right\Vert -\left\Vert \textrm{grad }L\right\Vert \right|d\mathcal{L}^{n}\\
 & \geq & \int_{B(0,R)}\left\Vert \textrm{grad }L\right\Vert d\mathcal{L}^{n}-\int\left\Vert Df\right\Vert d\mathcal{L}^{n}\end{eqnarray*}
so, by Lemma (\ref{lem:lemma2})\begin{eqnarray*}
\left\Vert \textrm{grad }L\right\Vert R^{n} & \leq & \int_{B(0,R)}\left\Vert Df\right\Vert d\mathcal{L}^{n}+\int_{B(0,R)}\left\Vert D(f+L)\right\Vert d\mathcal{L}^{n}\\
 & \leq & \sqrt{E}R^{n}+\left\Vert H_{\#}(T)\right\Vert \\
 & \leq & \sqrt{E}R^{n}+(\alpha+1)R^{n}.\end{eqnarray*}
Thus, by Proposition (\ref{lem:Main Lemma}), which implies that there
is a linear map $h$ for which $Exc(h_{\#}(T),\beta R)\leq c\beta^{2}Exc(T,R)$,
so that, when $k$ is given by $k:=l-L$, $k$ satisfies the excess
conditions (\ref{eq:excess bound}) and (\ref{eq:iterative excess bound})
of this Corollary and \[
\left\Vert \textrm{grad }k\right\Vert \leq\left\Vert \textrm{grad }l\right\Vert +\left\Vert \textrm{grad }L\right\Vert ,\]
which, for $\alpha>0$ sufficiently small satisfies the gradient bound
condition (\ref{eq:gradient bound}).
\hfill{}$\square$\end{pf}
The primary use of Proposition (\ref{lem:Main Lemma}) and its corollary
is in the following Lemma. Let $Exc(T,a,r)$ be the excess of $T$
over $B(a,r)$. 

\begin{lem}
\label{lem:lemma 17} There is a positive constant $E_{0}$ with the
following property. If $T$ is as in Proposition (\ref{lem:Main Lemma})
and \[
R+Exc(T,0,R)\leq E_{0},\]
then, for all $a,\, r$ with $\left|a\right|<R/2$, $r\leq R/2$ there
is a linear map $h(x,y)=(x,y-l(x))$ so that \[
Exc(h_{\#}(T),a,r)\leq C\left(\frac{r}{R}\right)^{2}Exc(T,0,R).\]
Moreover, $\left|grad(l)\right|\leq1/E_{0}$. 
\end{lem}
\begin{pf}
Since\[
Exc(T,a,R/2)\leq2^{n}Exc(T,0,R),\]
by replacing $E_{0}$by $E_{0}/2^{n}$ we see that it is sufficient
to prove the Lemma in the special case $a=0$. To do so, we apply
Proposition (\ref{lem:Main Lemma}) and Corollary (\ref{cor:Main-Lemma-tilted})
several times. Each time, we may need to use a smaller $\alpha$,
but since our iteration is finite there will be a sufficiently small
$\alpha$ to work for all steps simultaneously. We get linear maps
$h_{1},\, h_{2},\,\ldots, h_{s}$, $h_{i}(x,y)=(x,y-l_{i}(x))$ so
that \[
Exc((h_{i})_{\#}(T),0,\beta^{i}R)\leq c\beta^{2}Exc((h_{i-1})_{\#}(T),0,\beta^{i-1}R)\]
for $i=1,\ldots,s$, and \[
\left|grad(l_{i}-l_{i-1})\right|\leq\alpha^{-1}\sqrt{Exc((h_{i-1})_{\#}(T),0,\beta^{i-1}R)},\]
with $l_{0}=0$, and $\alpha$ satisfying the conditions of Proposition
(\ref{lem:Main Lemma}) or Corollary (\ref{cor:Main-Lemma-tilted})
for $(h_{i})_{\#}(T)$. Iterating the first inequality $s$ times,
\begin{eqnarray*}
Exc((h_{s})_{\#}(T),0,\beta^{s}R) & \leq & c\beta^{2}Exc((h_{s-1})_{\#}(T),0,\beta^{s-1}R)\\
 & \leq & c^{2}\beta^{4}Exc((h_{s-2})_{\#}(T),0,\beta^{s-2}R)\\
 & \vdots\\
 & \leq & c^{s}\beta^{2s}Exc(T,0,R).\end{eqnarray*}
So, choosing $s$ so that $\beta^{s+1}R<r\leq\beta^{s}R$, and $h=h_{s}$
we have the second claim of the Lemma. The second inequality from
Proposition (\ref{lem:Main Lemma}) or its corollary becomes \[
\left|grad(l_{i}-l_{i-1})\right|\leq\alpha^{-1}c^{(i-1)/2}\beta^{i-1}\sqrt{Exc(T,0,R)},\]
so\begin{eqnarray*}
\left|grad(l_{s})\right| & \leq\\
\sum_{i=1}^{s}\left|grad(l_{i}-l_{i-1})\right| & \leq & \sum_{i=1}^{s}\alpha^{-1}c^{(i-1)/2}\beta^{i-1}\sqrt{Exc(T,0,R)},\end{eqnarray*}
which, assuming at no loss in generality that $c\beta\sqrt{E}<1/2$,
is less than $2/\alpha$. Choosing $E_{0}=\alpha/2$ completes the
proof.
\hfill{}$\square$\end{pf}
\begin{prop}
\label{pro:lemma 18} With the hypotheses of Lemma (\ref{lem:lemma 17}),
the BV-carrier function $f(x)$ of $T$ is of class $C^{1}$ in $B(0,R/2)$. 
\end{prop}
\begin{pf}
Recall that $T\rest y_{i}=B(0,R)\rest f_{i}$. If $h(x,y)=(x,y-l(x)),$with
$l$ linear, then the corresponding function for $h_{\#}(T)$ is of
course $f-l$. Apply Lemma (\ref{lem:lemma2}) to the current $h_{\#}(T)\rest C(a,r)$,
implying from Lemma (\ref{lem:lemma 17}) that \[
\int_{B(a,r)}\left|Df-l\right|d\mathcal{L}^{n}\ll r^{n}\sqrt{C}\left(\frac{r}{R}\right)\sqrt{E},\]
for all $a\in B(0,R/2)$ and all $r\leq R/2$, for $l=l_{r,a}=D(l_{r,a})$,
where, from Lemma (\ref{lem:lemma 17}), note that the linear map
$l$ of that lemma, here denoted $l_{r,a}$, depends on the center
$a$ and radius $r$ of the ball. We need to show that the limit \[
l_{a}=\lim_{r\rightarrow0}l_{r,a}\]
exists for all $a\in B(0,R/2)$. Since \begin{eqnarray*}
\int_{B(a,r/2)}\left|l_{r,a}-l_{r/2,a}\right|d\mathcal{L}^{n} & \leq & \int_{B(a,r/2)}\left|Df-l_{r/2,a}\right|d\mathcal{L}^{n}+\int_{B(a,r)}\left|Df-l_{r,a}\right|d\mathcal{L}^{n}\\
 & \ll & r^{n}\sqrt{C}\left(\frac{r}{R}\right)\sqrt{E},\end{eqnarray*}
so, by the fact that the first integrand above is constant,\[
\left|l_{r,a}-l_{r/2,a}\right|\ll\sqrt{C}\left(\frac{r}{R}\right)\sqrt{E}.\]
Iterating that inequality and adding, \[
\left|l_{r,a}-l_{r/2^{n},a}\right|\ll\sqrt{C}\left(\frac{r}{R}\right)\sqrt{E}\sum_{j=0}^{\infty}2^{-j},\]
so, by the triangle inequality, $\left\{ l_{r/2^{n},a}\right\} $
is Cauchy in $\textrm{Hom}(\R^{n},\R^{j})$. Set $\hat{l}_{a}:=\lim l_{r/2^{n},a}$,
then $\left|l_{r,a}-\hat{l}_{a}\right|\ll r/R$, and so \[
l_{a}:=\lim_{r\rightarrow0}l_{r,a}=\hat{l}_{a}\]
exists. 

By a similar argument to the above, for $a,b\in B(0,R/2)$, with $\left|a-b\right|<r/2$,
\[
\left|l_{r,a}-l_{r,b}\right|\ll r/R,\]
and so \[
\left|l_{a}-l_{b}\right|\ll r/R\]
 if $\left|a-b\right|<r/2$, and so $a\mapsto l_{a}$ is continuous
in $B(0,R/2)$. 

From this it follows that $r^{-n}\int_{B(a,r)}\left|Df\right|d\mathcal{L}^{n}$
is uniformly bounded for $a\in B(0,R/2)$, $r<R/2$, and so the measure
$\left|Df\right|d\mathcal{L}^{n}$ is absolutely continuous, so that
we can write \[
Df\, d\mathcal{L}^{n}=\phi d\mathcal{L}^{n}\]
 for some $\phi\in L^{1}(B(0,R/2))$. Since $\phi\in L^{1}$, almost
every point in $B(0,R/2)$ is a Lebesgue point, and \[
\lim_{r\rightarrow0}r^{-n}\int_{B(a,r)}\left|\phi(x)-\phi(a)\right|d\mathcal{L}^{n}=0.\]
But we already have that \[
\lim_{r\rightarrow0}r^{-n}\int_{B(a,r)}\left|\phi(x)-l_{a}\right|d\mathcal{L}^{n}=0,\]
so that $\phi(a)=l_{a}$ almost-everywhere in $B(0,R/2)$, so $Df=\phi d\mathcal{L}^{n}$
with now $\phi$ continuous on $B(0,R/2)$, and so $f\in C^{1}(B(0,R/2))$.
\hfill{}$\square$\end{pf}
Similarly to the discussion in \cite[p. 129, lines 12-24]{Bombieri},
we have: 

Let $z\in Supp(T)$ be a point with an approximate tangent plane $Tan(Supp(T),z)$.
By the rectifiability theorem for currents and our lower bound on
density Proposition (\ref{pro:monotonicity}), we have that $Tan^{n}(\left\Vert T\right\Vert ,z)=Tan(Supp(T),z)$,
and is an $n$-dimensional vector space. If the oriented tangent plane
$\overrightarrow{T_{z}}$ is not vertical, that is, $\pi_{*}\overrightarrow{T_{z}}=\mathbf{e}$,
then there is a linear map $H(x,y):=(x,y-L(x))$ for which $\overrightarrow{H_{\#}(T)_{H(z)}}=\mathbf{e}_{0}$.
By Corollary (\ref{cor:Main-Lemma-tilted}), Lemma (\ref{lem:lemma 17})
and so Proposition (\ref{pro:lemma 18}) will apply to $H_{\# }(T)$
as well. By the monotonicity result, the density $\Theta(T,z)=1$
at each point. As a final assumption, assume that the tangent cone
$Tan(Supp(T),z)\subset V$, where $V\subset\R^{n+j}$ is an $n$-dimensional
plane. Since $Tan(Supp(T),Z)\supset Supp(\overrightarrow{T_{z}})$,
the plane $V=Supp(\overrightarrow{T_{z}})$ is not vertical. Apply
a shear-type linear map $H(x,y):=(x,y+L(x))$ so that $H(Tan(Supp(T),Z)=\R^{n}\times\{0\}\subset\R^{n+j}$.
Presume that coordinates are chosen so that $z=(0,0)$.

\begin{prop}
Under the conditions of Lemma (\ref{pro:lemma 18}), $Supp(T)$ is
a $C^{1}$, $n$-dimensional graph over some ball $B(0,r)$.
\end{prop}
\begin{pf}
It of course suffices to show that $Supp(H_{\#}(T))$ is a $C^{1}$
graph over $B(0,r)$. Since the tangent plane of $H_{\#}(T)$ over
$0$ is the horizontal plane in the coordinate system of the last
line of the previous paragraph, given $\eta>0$, there is an $r=r_{\eta}>0$
so that \[
Supp(H_{\#}(T)\rest C(0,r))\subset\left\{ \left|x\right|\leq r,\left|y\right|\leq\eta r\right\} =B(0,r)\times B(0,\eta r),\,\textrm{if }r\leq r_{\eta},\]
and \[
Supp(\partial(H_{\#}(T)\rest B(0,r)\times B(0,\eta r)))\subset\partial B(0,r)\times B(0,\eta r).\]
Once we show that \[
\lim_{r\rightarrow0}Exc(H_{\#}(T),0,r)=0,\]
then we can apply Lemma (\ref{lem:lemma 17}) and Proposition (\ref{pro:lemma 18})
to complete the proof of the present Proposition. 

Let $T$ be energy-minimizing among rectifiable sections, $T_{0}$
the current $B(0,r)\times\left\{ 0\right\} \subset B(0,r)\times B(0,\eta r)$,
and $F_{r}$ be the {}``fence'' obtained by connecting each element
of $(x,y)\in Supp(\partial H_{\#}(T)\rest B(0,r)\times B(0,\eta r))$
to $(x,0)\in Supp(\partial T_{0})$. Note that $T_{0}$ and $H_{\#}(T)\rest B(0,r)\times B(0,\eta r)-F_{r}$
have the same boundary, $\partial T_{0}$. It is easy to see (\cite[p. 363]{GMT},
or \cite[p. 128]{Bombieri}) that, \[
\left\Vert F_{r}\right\Vert \leq\left(\sup_{\partial T}\left|y\right|\right)\left\Vert \partial H_{\#}(T)\rest B(0,r)\times B(0,\eta r)\right\Vert ,\]
and, by slicing and the monotonicity formula, for a generic $\rho$,
$r<\rho<2r$ ($r<R/2$), there is a $C$ so that $\left\Vert \partial H_{\#}(T)\rest B(0,\rho)\times B(0,\eta\rho)\right\Vert \leq C\rho^{n-1}$.
Combining these two inequalities together, \[
\left\Vert F_{\rho}\right\Vert \leq C\eta\rho^{n}.\]

Since each penalty functional satisfies the ellipticity bounds (equation
(\ref{eq:ellipticity})),\[
\left[\left\Vert H_{\#}(T)\rest B(0,\rho)\times B(0,\eta r)-F_{\rho}\right\Vert -\left\Vert T_{0}\right\Vert \right]\leq\mathcal{F}_{\epsilon}(H_{\#}(T)\rest B(0,\rho)\times B(0,\eta\rho)-F_{\rho})-\mathcal{F}_{\epsilon}(T_{0}),\]
then so will the limiting functional $\mathcal{F}$. Then, by subadditivity,
and minimality of $T$, \begin{eqnarray*}
\left[\left\Vert H_{\#}(T)\rest B(0,\rho)\times B(0,\eta\rho)-F_{\rho}\right\Vert -\left\Vert T_{0}\right\Vert \right] & \leq & \mathcal{F}(H_{\#}(T)\rest B(0,\rho)\times B(0,\eta\rho)-F_{\rho})-\mathcal{F}(T_{0})\\
 & \leq & \left\Vert H\right\Vert ^{n}\left[\mathcal{F}(T\rest H^{-1}(B(0,\rho)\times B(0,\eta\rho)))\right.\\
 &  & \left.-\mathcal{F}(H_{\#}^{-1}(T_{0}+F_{\rho}))+2\mathcal{F}(H_{\#}^{-1}(F_{\rho}))\right]\\
 & \leq & 2\left\Vert H\right\Vert ^{n}\mathcal{F}(H_{\#}^{-1}(F_{\rho}))\\
 & \leq & 2C\eta\rho^{n}\end{eqnarray*}
and so \begin{eqnarray*}
Exc(H_{\#}(T),r) & = & \left.\left(\left\Vert H_{\#}(T)\rest B(0,r)\times B(0,\eta r)\right\Vert -\left\Vert T_{0}\right\Vert \right)\right/r^{n}\\
 & \leq & 2^{n}\left.\left(\left\Vert H_{\#}(T)\rest B(0,\rho)\times B(0,\eta\rho)\right\Vert -\left\Vert T_{0}\right\Vert \right)\right/\rho^{n}\\
 & \leq & 2^{n}\left.\left[\left\Vert F_{\rho}\right\Vert +\left\Vert T\rest B(0,\rho)\times B(0,\eta\rho)-F_{\rho}\right\Vert -\left\Vert T_{0}\right\Vert \right]\right/\rho^{n}\\
 & \leq & 2^{n}3C\eta.\end{eqnarray*}
Since, for any $\eta>0$ there is an $r>0$ sufficiently small so
that the conditions of Proposition (\ref{pro:lemma 18}) hold, the
conclusion of the Proposition holds.
\hfill{}$\square$\end{pf}
This proposition shows that the set of {}``good'' points in the
base manifold $M$, the set of points where there is a non-vertical
tangent space, is an open set, and on that open set the graph is of
class $C^{1}$. The next result completes the proof of the main theorem,
Theorem (\ref{thm:main}).

\begin{prop}
Let $T$ be an $n$-dimensional, mass-minimizing rectifiable section
in $\widetilde{\Gamma}(B)$ which is the limit of a sequence of penalty-minimizers.
Then, the projection $\pi(S)=Z$ onto $\mathbf{e}$ of the set $S$
of all points $y\in Supp(T)$ so that the oriented tangent cone is
not a plane, or where $T(T,y)$ has a vertical direction a closed
set of Hausdorff $n$-dimensional measure 0 in $M$.
\end{prop}
\begin{pf}
The previous section shows that the set of points with non-vertical
tangent planes is open in $T$, and projects to an open set. So, the
set of points with no tangent plane, or with one having vertical directions,
is closed. The set of points with no tangent plane is of measure $0$
in any countably-rectifiable integer-multiplicity current, by \cite[3.2.19]{GMT}.
Assume there is a set $S_{0}\subset S$ of points of $T$ with $Z_{0}:=\pi(S_{0})$
of positive Hausdorff $n$-dimensional measure, and with vertical
tangent planes. For all but a set of Hausdorff $n$-dimensional measure
0 in $S_{0}$, the density of the set will be 1 as well \cite[3.2.19]{GMT}.
For any such $z\in S_{0}$, given any $\epsilon>0$, there is a $\delta_{z,\epsilon}>0$
so that the ratio of the measure of of the projection onto $M$ of
$S_{0}\cap B(z,\delta_{z,\epsilon})$ to that of the ball of radius
$\delta_{z,\epsilon}$centered at $\pi(z)$ in $M$ will be less than
$\epsilon$,\[
\epsilon>\frac{\mathcal{H}^{n}(\pi(S_{0}\cap B_{B}(z,\delta_{z,\epsilon})))}{\mathcal{H}^{n}(B_{M }(\pi(z),\delta_{z,\epsilon}))}>\frac{\mathcal{H}^{n}(\pi(S_{0}\cap B_{B}(z,\delta_{z,\epsilon})))}{\omega_{n}\delta_{z,\epsilon}^{n}/2}\]
 since the tangent plane is vertical. But the measure $\mathcal{H}^{n}(S_{0}\cap B_{B}(z,\delta_{z,\epsilon}))>\omega_{n}\delta_{z,\epsilon}^{n}/2$
for small enough $\delta_{z,\epsilon}$ since the density is one,
so by the Besicovitch covering theorem \[
\mathcal{H}^{n}(S_{0})>\frac{1}{\epsilon}\mathcal{H}^{n}(Z_{0}).\]
Since this must be true for any $\epsilon>0$, it would contradict
the fact that $T$ has finite mass if $\mathcal{H}^{n}(Z_{0})>0$.
\hfill{}$\square$\end{pf}

\end{article}

\begin{thebibliography}{99}
\bibitem{Almgren}F. J. Almgren, Jr. \textit{Q-valued functions minimizing Dirichlet's
integral and the regularity of area minimizing rectifiable currents
up to codimension two ,} Bull. (New Series) Am. Math. Soc., \textbf{8}
(1983), 327--328. 
\bibitem{Bombieri}E. Bombieri, \emph{Regularity theory for almost minimal currents},
Arch. Rat. Mech. Anal. \textbf{78} (1982), 99-130.
\bibitem{Coventry}A. Coventry, Research reports (Mathematics), number CMA-MRR 45- 98,
Australian National University Publications, 1998.
\bibitem{Evans and Gariepy}L.C. Evans and R. F. Gariepy, \emph{Measure Theory and Fine Properties
of Functions}. Studies in Advanced Mathematics. CRC Press, 1992.
\bibitem{GMT}Herbert Federer, \textit{Geometric Measure Theory,} Springer-Verlag
1969. 
\bibitem{j}D. L. Johnson, \textit{K\"{a}hler submersions and holomorphic connections,}
Jour. Diff. Geo. \textbf{15} (1980) , 71--79.
\bibitem{js-1}D. L. Johnson and P. Smith, \textit{Regularity of volume-minimizing
graphs,} Indiana University Mathematics Journal, \textbf{44} (1995),
45--85. 
\bibitem{dim-3}D. L. Johnson and P. Smith, \textit{Regularity of mass-minimizing
one-dimensional foliations,} Analysis and Geometry on Foliated Manifolds,
Proceedings of the VII International Colloquium on Differential Geometry,
(1994), World Scientific, 81--98. 
\bibitem{Lin}Lin Fanghua and Yang Xiaoping, \emph{Geometric Measure Theory, an
Introduction}, International Press, Boston, 2002.
\bibitem{MS}J. D. Moore and R. Schlafly, \textit{On equivariant isometric embeddings,}
Mathematische Zeitschrift \textbf{173} (1980), 119--133. 
\bibitem{Morgan}F. Morgan, \emph{Geometric Measure Theory, A Beginner's Guide}, Academic
Press, 1988; second edition, 1995.
\bibitem{Morrey}Morrey, C. B., \emph{Multiple Integrals in the Calculus of Variations},
Springer-Verlag, 1966.
\bibitem{Roubicek}T. Roub\'{i}\v{c}ek, \emph{Relaxation in Optimization Theory and
Variational Calculus}, Walter de Gruyter, 1999.
\bibitem{Sasaki}Takeshi Sasaki, \textit{On the differential geometry of tangent bundles
of Riemannian manifolds}, T\^ ohoku Math. J. \textbf{10} (1958),
338--354.
\bibitem{Serrin}James Serrin, \textit{On the definition and properties of certain
variational integrals,} Trans. Amer. Math. Soc. \textbf{101} (1961),
139--167. 
\bibitem{Struwe}Michael Struwe, \emph{Variational methods : applications to nonlinear
partial differential equations and Hamiltonian systems,} Springer-Verlag
Ergebnisse der Mathematik und ihrer Grenzgebiete v. 34, 2000.
\end{thebibliography}
\end{document}